\documentclass[11pt]{article}
\usepackage{amsmath}
\usepackage{amssymb}
\usepackage{amsthm}
\usepackage{amscd}
\usepackage{amsfonts}
\usepackage{dsfont}
\usepackage{graphicx}
\usepackage{subcaption}
\usepackage{fancyhdr}
\usepackage{mathtools}
\usepackage{cite}
\usepackage[dvipsnames]{xcolor}
\usepackage[ruled]{algorithm2e}
\usepackage{ifthen}
\usepackage{algorithmic}
\newtheorem{theorem}{Theorem}
\newtheorem{lemma}{Lemma}

\newtheorem{proposition}{Proposition}
\newtheorem{corollary}{Corollary}

\newtheorem{definition}{Definition}
\newtheorem{assumption}{Assumption}
\DeclareMathOperator*{\argmin}{arg\,min}
\DeclareMathOperator*{\argmax}{arg\,max}
\usepackage{empheq}
\def\scr#1{{\cal #1}} 
\def\eq#1{\begin{equation}#1\end{equation}}
\newcommand{\R}{{\rm I\!R}}
\newcommand{\bbb}{\mathbb}

\def\rep#1{(\ref{#1})}
\newcommand{\1}{\mathbf{1}}

\newcommand{\dfb}{\stackrel{\Delta}{=}}
\usepackage{hyperref}

\def\##1\#{\begin{align}#1\end{align}}
\def\$#1\${\begin{align*}#1\end{align*}}

\def\qed{ \rule{.08in}{.08in}}

\usepackage{parskip}
\title{Resilient Distributed Optimization 
}

\author{Jingxuan Zhu \hspace{.3in} Yixuan Lin \hspace{.3in} Alvaro Velasquez \hspace{.3in} Ji Liu\thanks{
J.~Zhu and Y.~Lin are with the Department of Applied Mathematics and Statistics at Stony Brook University (\texttt{\{jingxuan.zhu,yixuan.lin.1\}@stonybrook.edu}).
A.~Velasquez is with the Department of Computer Science at University of Colorado Boulder (\texttt{alvaro.velasquez@colorado.edu}).
J.~Liu is with the Department of Electrical and Computer Engineering at Stony Brook University
(\texttt{ji.liu@stonybrook.edu}).
}
}

\begin{document}
\date{}

\maketitle
\thispagestyle{empty}

\begin{abstract}
This paper considers a distributed optimization problem in the presence of Byzantine agents capable of introducing
untrustworthy information into the communication network. A resilient distributed subgradient algorithm is proposed based on graph redundancy and objective redundancy. It is shown that the algorithm causes all non-Byzantine agents' states to asymptotically converge to the same optimal point under appropriate assumptions. A partial convergence rate result is also provided. 
\end{abstract}

\section{Introduction}

Distributed optimization has attracted considerable attention and achieved remarkable success in both theory and practice. The distributed convex optimization problem was first studied in \cite{nedic2009distributed} where a distributed subgradient algorithm was proposed. After this, various distributed optimization algorithms have been crafted and studied; see survey papers \cite{yang2019survey,nedic2018distributed,molzahn2017survey}. Distributed optimization techniques are also widely applied to decentralized deep learning \cite{LianNIPS2017}.

Information exchange between neighboring agents is necessary for a multi-agent network for distributed optimization. However, agents' states may be corrupted and they may not adhere to the designed algorithm due to faulty processes or external attacks. An agent is called {\em Byzantine} if it updates its state in an arbitrary, unknown manner, and can send conflicting values to different neighbors \cite{byzantine}. Such attacking agents can know global information of the network, play arbitrarily and strategically, and even be coordinated. 
Consider a network of agents in which Byzantine agents exist. An ideal resilient algorithm is the one which can lead non-Byzantine (or normal) agents to cooperatively solve the
corresponding distributed optimization problem in the presence of Byzantine agents as if they do not exist. Such a resilient algorithm is highly desirable for the safety and security of multi-agent systems as faulty processes and external attacks are inevitable.

Resilient distributed optimization has recently received increasing attention, probably originating from the work of \cite{sulili}.
Almost all the existing works cannot guarantee full resilience; what they can guarantee is all normal agents' states converge to a bounded neighborhood of the desired optimal point whose bound is not controllable \cite{kuwaranancharoen2020byzantine,8759887,fang2022bridge,clipping,tianyi}, or an optimal point of an unspecified convex combination of all normal agents' objective functions \cite{sulili,sundaram2018distributed,su2020byzantine}, or a convex combination of all normal agents' local optimal points \cite{chengcheng}. 
The only exceptions are \cite{gupta2021byzantine,gupta2021byzantine_3,ElioUsai}  in which the underlying communication graph is assumed to be a complete graph, namely, each agent is allowed to communicate with all other agents. All \cite{gupta2021byzantine,gupta2021byzantine_3,ElioUsai} rely on the idea of ``objective function redundancy''. The idea has also been applied to the federated setting and achieved full resilience \cite{liu2021approximate,gupta2021byzantine_2}. In the federated setting, a central coordinator agent is able to communicate with all worker agents, which is more or less equivalent to a complete graph in the distributed setting (or sometimes called decentralized setting). It is worth noting that \cite{sulili,sundaram2018distributed,su2020byzantine,chengcheng,ElioUsai} only consider special one-dimensional optimization.

Resilient distributed optimization is also related to resilient federated optimization/learning in the coordinator-workers setting (e.g., \cite{gupta2021byzantine_2,9248056,9030051}), which has attracted increasing attention recently. The key problem is how the central coordinator aggregates the received information to eliminate or attenuate the effects of Byzantine worker agents. Various Byzantine-resilient information aggregation methods have been proposed for high-dimensional optimization/learning, focusing on stochastic gradient descent (SGD). Notable among them are \cite{csx,blanchard,Alistarh,suxu,chen2018draco}, just to name a few; see an overview of recent developments in this area in \cite{9084329}. It is doubtable that these methods can be applied to achieve full resilience in the distributed setting.

From the preceding discussion, and to the best of our knowledge, a fully resilient distributed optimization algorithm for general non-complete communication graphs does not exist, even for one-dimensional optimization problems. This gap is
precisely what we study in this paper. We consider a distributed convex optimization problem in the presence of Byzantine agents and propose a fully resilient distributed subgradient algorithm based on the ideas of objective redundancy (cf. Definition \ref{def:redundant}) and graph redundancy (cf. Definition \ref{def:resilient}). The algorithm is shown to cause all non-Byzantine agents' states to asymptotically converge to the same desired optimal point under appropriate assumptions. The proposed algorithm works theoretically for multi-dimensional optimization but practically not for  high-dimensional optimization, as will be explained and discussed in the concluding remarks.

This work is motivated by two recent ideas. The first is the quantified notion of objective function redundancy proposed in \cite{reviewer} where a couple of different definitions of objective redundancy are studied, based on which fully resilient distributed optimization algorithms have been crafted either for a federated setting \cite{reviewer,liu2021approximate,gupta2021byzantine_2} or a distributed setting over complete graphs \cite{gupta2021byzantine,gupta2021byzantine_3,ElioUsai}; such redundancy has been shown necessary for achieving full resilience in multi-agent optimization \cite{liu2021approximate}. 
We borrow one notation in \cite{reviewer} and further develop it. It is worth emphasizing that the results in \cite{gupta2021byzantine,gupta2021byzantine_3,ElioUsai} rely on objective redundancy among non-Byzantine agents, whereas ours depend on objective redundancy among all agents. This subtle difference is important for equipping a multi-agent network with a certain level of redundancy at a network design stage  as which agents are non-Byzantine cannot be assumed a priori.

The second idea is so-called ``Byzantine vector consensus'' \cite{Va13,Va14} whose goal is, given a set of both Byzantine and non-Byzantine vectors, to pick a vector lying in the convex hull of the non-Byzantine vectors, based on Tverberg's theorem \cite{tverberg,tverbergpoint}.
The idea has been very recently improved in \cite{resilientconstrained} which can be used to achieve resilient multi-dimensional consensus exponentially fast. Exponential consensus is critical in the presence of diminishing disturbance \cite{cdc16consensus}. We are prompted by this improved idea and make use of a resilient vector picking process, simplified from that of \cite[Algorithm 1]{resilientconstrained}. There are other recent approaches appealing to the idea of centerpoint \cite{centerpoint,yan2022resilient}. We expect that these approaches can also be applied to resilient optimization, provided that exponential consensus is guaranteed, e.g., in \cite{yan2022resilient}.

\section{Problem Formulation}

Consider a network consisting of $n$ agents, labeled $1$ through $n$ for the purpose of presentation. The agents are not aware of such global labeling, but can differentiate between their neighbors. The neighbor relations among the $n$ agents are characterized by a directed graph $\bbb{G} = (\mathcal{V},\mathcal{E})$ whose
vertices correspond to agents and whose directed edges (or arcs) depict neighbor relations, where $\mathcal{V}=\{1,\ldots,n\}$ is the vertex set and $\mathcal{E}\subset\mathcal{V} \times \mathcal{V}$ is the directed edge set.\footnote{We use $\scr{A} \subset \scr{B}$ to denote that $\scr A$ is a subset of $\scr B$.}
Specifically, agent $j$ is a neighbor of agent $i$ if $(j,i)\in\scr{E}$. 
Each agent can receive information from its neighbors. Thus, the directions of edges represent the
directions of information flow. 
We use $\mathcal{N}_i$ to denote the neighbor set of agent $i$, excluding $i$, i.e.,
$\mathcal{N}_i = \{ j \in \mathcal{V}  : ( j, i ) \in \mathcal{E} \}$.


Each agent $i\in\scr V$ has a ``private'' convex (not necessarily differentiable) objective function, $f_i:\R^d\rightarrow\R$, only known to agent $i$. There exist Byzantine agents in the network which are able to transmit arbitrary values to others and capable of sending conflicting values to different neighbors at any time. The set of Byzantine agents is denoted by $\scr F$ and the set of normal (non-Byzantine) agents is denoted by $\scr H$. Which agents are Byzantine is unknown to normal agents. 
It is assumed that each agent may have at most $\beta$ Byzantine neighbors.

The goal of the normal agents is to cooperatively minimize the objective functions
$$f_{\scr H}(x)=\sum_{i\in\scr H} f_i(x) \;\;\;\;\; {\rm and} \;\;\;\;\;f(x)=\sum_{i\in\scr V} f_i(x).$$
We will show that minimizing the above two objective functions can be achieved simultaneously with appropriate redundancy in objective functions (cf. Definition \ref{def:redundant} and Corollary \ref{coro:optset}). 
It is assumed that the set of optimal solutions to $f$, denoted by $\scr X^*$, is nonempty and bounded.

Since each $f_i$ is not necessarily differentiable, the gradient descent method may not be applicable. Instead, the subgradient method \cite{subgradient} can be applied. For a convex function $h : \R^d\rightarrow \R$, a vector $g\in\R^d$ is called a subgradient of $h$ at point $x$ if
\eq{
h(y)\ge h(x) + g^\top (y-x) \;\; {\rm for \; all} \;\; y\in\R^d.
\label{eq:subgradient}}
Such a vector $g$ always exists and may not be unique. In the case when $h$ is differentiable at point $x$, the subgradient $g$ is unique and equals $\nabla h(x)$, the gradient of $h$ at $x$. Thus, the subgradient can be viewed as a generalization of
the notion of the gradient. 
From \eqref{eq:subgradient} and the Cauchy-Schwarz inequality, 
$h(y) - h(x) \ge - G \| y-x\|$,
where $G$ is an upper bound for the 2-norm of the subgradients of $h$ at both $x$ and $y$. We will use this fact without special mention in the sequel.
Throughout this paper, we use $\|\cdot\|$ for the 2-norm.


The subgradient method was first proposed in \cite{subgradient} and the first distributed subgradient method was proposed in \cite{nedic2009distributed} for undirected graphs. Its extension to directed graphs has been studied in  \cite{nedic} and recently further analyzed in \cite{cdc22push}.

\subsection{Redundancy}

To make the resilient distributed optimization problem solvable,
certain redundancy is necessary. We begin with objective redundancy.

\begin{definition}\label{def:redundant}
An $n$-agent network is called $k$-redundant, $k\in\{0,1,\ldots,n-1\}$, if for any subsets $\scr S_1,\scr S_2\subset\scr V$ with $|\scr S_1|= |\scr S_2| = n - k$, there holds\footnote{We use $|\scr S|$ to denote the cardinality of a set $\scr S$.} 
\begin{align*}
    \argmin_{x}\sum_{i\in\scr S_1}
    f_i(x)=\argmin_{x}\sum_{i\in\scr S_2}f_i(x).
\end{align*} 
\end{definition}


The above definition of objective redundancy originated in \cite[Definition 2]{reviewer}. It has the following properties.

\begin{lemma}\label{le:equaloptimalset}
If an $n$-agent network is $k$-redundant, then for any subsets $\scr S,\scr L\subset\scr V$ with $|\scr S|= n-k$ and $|\scr L|\ge n-k$, there holds
\[\argmin_{x}\sum_{i\in\scr S}
    f_i(x)=\argmin_{x}\sum_{i\in\scr L}f_i(x).\]
\end{lemma}

{\bf Proof of Lemma \ref{le:equaloptimalset}:}
Let $\scr Z= \argmin_{x}\sum_{i\in\scr S}
    f_i(x)$ and $\scr Q=\{\scr P\;:\;\scr P\subset \scr L,\; |\scr P|=n- k\}$.
From Definition~\ref{def:redundant}, $\argmin_{x}\sum_{i\in\scr P}f_i(x)=\scr Z$ for any $\scr P\in\scr Q$. 
For each $i\in\scr L$, let $\scr Q_i=\{\scr P\;:\;\scr P\subset \scr L,\; |\scr P|=n- k,\; i\in\scr P\}$. It is easy to see that for each $i\in\scr L$,\footnote{We use $\binom{n}{k}$ to  denote the number of $k$-combinations from a set of $n$ elements.}
$$|\scr Q_i| = q \dfb \binom{|\scr L|-1}{n- k-1}.$$
Then, 
    \begin{align}
        \sum_{\scr P\in\scr Q}\sum_{i\in\scr P}f_i(x) = q\sum_{i\in\scr L}f_i(x).\label{eq:qqq}
    \end{align}
Pick any $z\in\scr Z$. From \eqref{eq:qqq},     
\begin{align*}
    \min_x q\sum_{i\in\scr L}f_i( x)= \min_x\sum_{\scr P\in\scr Q}\sum_{i\in\scr P}f_i(x)
    \ge\sum_{\scr P\in \scr Q}\min_x\sum_{i\in\scr P}f_i(x)=\sum_{\scr P\in\scr Q}\sum_{i\in\scr P}f_i(z) = q\sum_{i\in\scr L}f_i(z),
\end{align*}
which implies that $z\in\argmin_x\sum_{i\in\scr L}f_i(x)$, and thus $\scr Z\subset\argmin_x\sum_{i\in\scr L}f_i(x)$. 

To prove the lemma, it is sufficient to prove that $\argmin_x\sum_{i\in\scr L}f_i(x)\subset \scr Z$. Suppose that, to the contrary, there exists a $y$ such that  $y\in\argmin_x\sum_{i\in\scr L}f_i(x)$ and $y \notin Z$. Since $y,z \in \argmin_x\sum_{i\in\scr L}f_i(x)$, from \eqref{eq:qqq},
\begin{align*}
    \sum_{i\in\scr L}f_i(y) = \sum_{i\in\scr L}f_i(z)= \frac{1}{q}\sum_{\scr P\in\scr Q}\sum_{i\in\scr P}f_i(z)< \frac{1}{q}\sum_{\scr P\in\scr Q}\sum_{i\in\scr P}f_i(y)=\sum_{i\in\scr L}f_i(y),
\end{align*}
which is impossible. Therefore, $\argmin_x\sum_{i\in\scr L}f_i(x)\subset \scr Z$.
\hfill$\qed$

The following corollaries are immediate consequences of Lemma~\ref{le:equaloptimalset}.

\begin{corollary}\label{coro:optset}
If an $n$-agent network is $k$-redundant, then for any subsets $\scr S\subset\scr V$ with $|\scr S|\ge n-k$, there holds
\[\argmin_{x}\sum_{i\in\scr S}
    f_i(x)=\scr X^*.\]
\end{corollary}

\begin{corollary}
If an $n$-agent network is $(k+1)$-redundant with $k\ge 0$, then it is $k$-redundant.
\end{corollary}

We also need redundancy in graph connectivity. 

A vertex $i$ in a directed graph $\bbb G$ is called a root of $\bbb G$ if for each other vertex $j$ of $\bbb G$, there is a directed path from $i$ to $j$. Thus, $i$ is a root of $\bbb G$ if it is the root of a directed spanning tree of $\bbb G$. We will say that $\bbb G$ is rooted at $i$ if $i$ is in fact a root. 
It is easy to see that a rooted graph $\bbb G$ has a unique strongly connected component whose vertices are all roots of $\bbb G$.

\begin{definition}\label{def:resilient}
An $(r,s)$-reduced graph of a directed graph $\bbb G$ with $n$ vertices, with $r,s\ge 0$ and $r+s\le n-1$, is a subgraph of $\bbb G$ obtained by first picking any vertex subset $\scr S\subset\scr V$ with $|\scr S|=n-r$ and then removing from each vertex of the subgraph induced by $\scr S$, $\bbb G_{\scr S}$, arbitrary $s$ incoming edges in $\bbb G_{\scr S}$. 
A directed graph $\bbb G$ is called $(r,s)$-resilient if all its $(r,s)$-reduced graphs are rooted.
\end{definition}

It is easy to see that if a directed graph is $(r_1,s_1)$-resilient, then for any nonnegative $r_2\le r_1$ and $s_2\le s_1$, the graph is also $(r_2,s_2)$-resilient.


In the case when $r=s=\beta$, the resilient graph is equivalent to rooted ``reduced graph'' in \cite{Va12} which was used to guarantee resilient one-dimension consensus; see Definition 4 and Theorem 2 in \cite{Va12}. Thus, the definition here can be viewed as a simple generalization of the rooted ``reduced graph''.

Definition \ref{def:resilient} implicitly requires that each vertex of an $(r,s)$-resilient graph have at least $r+s$ neighbors. More can be said.

\begin{lemma}\label{le:enoughneighbor}
If a directed graph is $(r,s)$-resilient, then each of its vertices has at least $(r+s+1)$ neighbors.
\end{lemma}


{\bf Proof of Lemma \ref{le:enoughneighbor}:}
Suppose that, to the contrary, there exists a vertex $i$ in $\bbb G$ whose $|\scr N_i|\le r+s$. If $|\scr N_i|< r+s$, it is easy to see that $\bbb G$ does not satisfy Definition \ref{def:resilient}. We thus consider the case when $|\scr N_i|= r+s$. Let $\scr R$ be the set of arbitrary $r$ neighbors of vertex $i$, and $\scr S=\scr V\setminus\scr R$, where $\scr V$ is the vertex set of $\bbb G$.\footnote{We use $\scr{A}\setminus \scr{B}$ to denote the set of elements that are in $\scr{A}$ but not in $\scr{B}$.} It is clear that $|\scr S|=n-r$, and in the subgraph induced by $\scr S$, $\bbb G_{\scr S}$, vertex $i$ has exactly $s$ neighbors. Then, after vertex $i$ removes $s$ incoming edges in $\bbb G_{\scr S}$, and each out-neighbor\footnote{A vertex $i$ is called an out-neighbor of vertex $j$ if the latter is a neighbor of the former.} of vertex $i$ in $\bbb G_{\scr S}$, if any, removes its incoming edge from $i$,
vertex $i$ becomes isolated. But it is impossible for an $(r,s)$-resilient graph.
\hfill$\qed$

\section{Algorithm}

To describe our algorithm, we need the following notation. 

Let $\scr A_i$ denote the collection of all those subsets of $\scr N_i$ whose cardinality is $(d+1)\beta + 1$. It is obvious that the number of all such subsets is
\eq{a_i \dfb \binom{|\scr N_i|}{(d+1)\beta + 1},\label{eq:a_i(t)}}
and label them $\scr A_{i1},\ldots, \scr A_{ia_i}$.
For each $j\in\{1,\ldots,a_i\}$, let $\scr B_{ij}$ denote the collection of all those subsets of $\scr A_{ij}$ whose cardinality is $d\beta + 1$.
For any subset of agents $\scr S\subset \scr V$, let $\scr C_{\scr S}(t)$ denote the convex hull of all $x_i(t)$, $i\in\scr S$.

{\bf Algorithm:} At each discrete time $t\in\{0,1,2,\ldots\}$, each agent $i$ first picks an arbitrary point 
\eq{y_{ij}(t)\in\bigcap_{\scr S\in\scr B_{ij}} \scr C_{\scr S}(t)\label{eq:hull}}
for each $j\in\{1,\ldots,a_i\}$, and 
then updates its state by setting
\begin{align}
    v_i(t) &= \frac{1}{1+a_i}\Big(x_i(t)+\sum_{j=1}^{a_i}y_{ij}(t)\Big),\label{eq:v_i(t)_ori}\\
    x_i(t+1) &= v_i(t)-\alpha(t)g_i(v_i(t)),\label{eq:x_ori}
\end{align}
where $\alpha(t)$ is the stepsize, and $g_i(\cdot)$ is a subgradient of $f_i(\cdot)$.
\hfill$\Box$

In the special one-dimensional case with $d=1$, it is not hard to check that the steps \eqref{eq:hull} and \eqref{eq:v_i(t)_ori} simplifies to the resilient scalar consensus algorithm in \cite{Va12}, which is essentially equivalent to the trimmed mean method and has been improved in \cite{LeZhKoSu13}.


The convergence and correctness of the proposed algorithm rely on the following assumptions.


\begin{assumption}\label{a:interior}
$\scr X^*$ has a nonempty interior.
\end{assumption}

It is easy to see that Assumption \ref{a:interior} implies that $f(x)$ is differentiable at any $x\in{\rm int}(\scr X^*)$, where ${\rm int}(\cdot)$ denotes the interior of a set. More can be said. 

\begin{lemma}\label{le:existenceofgradient}
Under Assumption \ref{a:interior}, if the $n$-agent network is $k$-redundant with $k\ge 1$, then $f_i(x)$ is differentiable at $x$ with $\nabla f_i(x)=0$ for all $i\in\scr V$ and $x\in{\rm int}(\scr X^*)$.
\end{lemma}


{\bf Proof of Lemma \ref{le:existenceofgradient}:}
Since ${\rm int}(\scr X^*)$ is nonempty, for any $x^*\in{\rm int}(\scr X^*)$, there exist a positive number $r$ and an open ball in ${\rm int}(\scr X^*)$ centered at $x^*$ with radius $r$, denoted as $\scr B(x^*,r)\subset {\rm int}(\scr X^*)$. Let $h_j$ be a vector in $\R^d$ whose $j$th entry is $\epsilon$ and the remaining entries all equal zero. Since $x^*+h_j\in\scr B(x_0,r)\subset {\rm int}(\scr X^*)$ for sufficiently small $\epsilon$,  
\begin{align}\label{eq:f_ilimit=0}
    \frac{\partial }{\partial x_j}f(x^*)=\lim_{\epsilon\rightarrow 0}\frac{f(x^*+h_j)-f(x^*)}{\epsilon}=
\lim_{\epsilon\rightarrow 0}\frac{\sum_{i\in\scr V}(f_i(x^*+h_j)-f_i(x^*))}{\epsilon}
=0.
\end{align}
For each $i\in\scr V$, since $f_i(x)$ is convex, 
both $\lim_{\epsilon\rightarrow 0^-}\frac{f_i(x^*+h_j)-f_i(x^*)}{\epsilon}$ and $\lim_{\epsilon\rightarrow 0^+}\frac{f_i(x^*+h_j)-f_i(x^*)}{\epsilon}$ exist and $\lim_{\epsilon\rightarrow 0^-}\frac{f_i(x^*+h_j)-f_i(x^*)}{\epsilon}\le\lim_{\epsilon\rightarrow 0^+}\frac{f_i(x^*+h_j)-f_i(x^*)}{\epsilon}$ for all $j\in\{1,\ldots,d\}$ \cite[Theorem 24.1]{Rockafellar+2015}.
It follows that 
\[\sum_{k\in\scr V}\lim_{\epsilon\rightarrow 0^-}\frac{f_k(x^*+h_j)-f_k(x^*)}{\epsilon}\le\sum_{k\in\scr V}\lim_{\epsilon\rightarrow 0^+}\frac{f_k(x^*+h_j)-f_k(x^*)}{\epsilon}.\]
Note that from \eqref{eq:f_ilimit=0},
\[\sum_{k\in\scr V}\lim_{\epsilon\rightarrow 0^-}\frac{f_k(x^*+h_j)-f_k(x^*)}{\epsilon}
=\sum_{k\in\scr V}\lim_{\epsilon\rightarrow 0^+}\frac{f_k(x^*+h_j)-f_k(x^*)}{\epsilon}.\]
Thus,
$\lim_{\epsilon\rightarrow 0^-}\frac{f_i(x^*+h_j)-f_i(x^*)}{\epsilon}=\lim_{\epsilon\rightarrow 0^+}\frac{f_i(x^*+h_j)-f_i(x^*)}{\epsilon}$, i.e., $\partial f_i(x^*)/\partial x_j$ exists for all $i\in\scr V$ and $j\in\{1,\ldots,d\}$.

To proceed, let $h_i(x)=\sum_{k\in\scr V,\;k\neq i}f_k(x)$ for all $i\in\scr V$. 
From Corollary \ref{coro:optset}, 
$\argmin_{x}h_i(x)=\scr X^*$.
Since $x^*\in{\rm int}(\scr X^*)$, 
both $f(x)$ and $h_i(x)$ are differentiable at $x^*$, implying that $\frac{\partial f}{\partial x_j}(x^*)=\frac{\partial h_i}{\partial x_j}(x^*)=0$ for all $i\in\scr V$ and $j\in\{1,\ldots,d\}$. Since $f_i(x)=f(x)-h_i(x)$, $\frac{\partial f_i}{\partial x_j}(x^*)=0$ for all $i\in\scr V$ and $j\in\{1,\ldots,d\}$. Note that this holds for all $x^*\in{\rm int}(\scr X^*)$. From \cite[Section 8.4.2]{zorich2004mathematical}, $f_i(x)$ is differentiable at $x^*$ with $\nabla f_i(x^*)=0$ for all $i\in\scr V$.
\hfill$\qed$


Lemma \ref{le:existenceofgradient} has the following important implication. 

\begin{corollary}\label{le:localoptimal}
Under Assumption~\ref{a:interior}, if the $n$-agent network is $k$-redundant with $k\ge 1$, then for all $i\in\scr V$, 
$$\scr X^*\subset \argmin_{x}f_i(x).$$
\end{corollary}

Corollary~\ref{le:localoptimal} immediately implies that 
\[\bigcap_{i\in\scr V}\argmin_{x}f_i(x)=\scr X^*.\]

{\bf Proof of Corollary \ref{le:localoptimal}:}
Suppose that, to the contrary, there exist $x^*\in\scr X^*$ and $i\in\scr V$ such that $x^*\notin\argmin_{x}f_{i}(x)$. Pick a $z\in {\rm int}(\scr X^*)$. From Lemma \ref{le:existenceofgradient},  $z\in\argmin_{x}f_i(x)$.
It is then clear that $f_{i}(x^*)>f_{i}(z)$.
Let $h_i(x)=\sum_{k\in\scr V,\;k\neq i}f_k(x)$. 
From Corollary \ref{coro:optset}, 
$\argmin_{x}h_i(x)=\scr X^*$, and thus $h_{i}(x^*)=h_{i}(z)$.
It follows that 
$f(x^*) = f_{i}(x^*)+h_{i}(x^*)>f_{i}(z)+h_{i}(z)=f(z)$,
which contradicts the fact that $x^*\in\scr X ^*$.
\hfill$\qed$

\begin{assumption}\label{assum:boundedsubgradient}
The subgradients of all $f_i$, $i\in\scr V$, are uniformly  bounded, i.e., there exists a positive number $D$ such that $\|g_i(x)\|\le D$ for all $i\in\scr V$ and $x\in\R^d$.
\end{assumption}

\begin{assumption} \label{assum:step-size}
    The step-size sequence $\{\alpha(t)\}$ is positive, non-increasing, and satisfies $\sum_{t=0}^\infty \alpha(t) = \infty$ and $\sum_{t=0}^\infty \alpha^2(t) < \infty$. 
\end{assumption}

The above two assumptions are standard for subgradient methods. 

To state our main results, we need the following concepts. For a directed graph $\bbb G$, we use $\scr R_{r,s}(\bbb G)$ to denote the set of all $(r,s)$-reduced graphs of $\bbb G$. For a rooted graph $\bbb G$, we use $\kappa(\bbb G)$ to denote the size of the unique strongly connected component whose vertices are all roots of $\bbb G$; in other words, $\kappa(\bbb G)$ equals the number of roots of $\bbb G$. For any $(r,s)$-resilient graph $\bbb G$, let 
$$\kappa_{r,s}(\bbb G)\dfb\min_{\bbb H\in\scr R_{r,s}(\bbb G)} \kappa(\bbb H).$$
which is well-defined and denotes the smallest possible number of roots in any $(r,s)$-reduced graphs of $\bbb G$. 


\begin{theorem}\label{thm:main}
Under Assumptions~\ref{a:interior}--\ref{assum:step-size}, if $\bbb G$ is $(\beta,d\beta)$-resilient and the $n$-agent network is $(n-\kappa_{\beta,d\beta}(\bbb G))$-redundant, then all $x_i(t)$, $i\in\scr H$ will asymptotically reach a consensus at a point in $\scr X^*$. 
\end{theorem}

The following example shows that $(n-\kappa_{\beta,d\beta}(\bbb G))$-redundancy is necessary. 
For simplicity, set $d=1$. Consider a 4-agent network whose neighbor graph is the 4-vertex complete graph $\bbb C$, which is $(1,1)$-resilient. Suppose that agent 4 is Byzantine and the other three are normal. It is possible that, with a carefully crafted attack strategy of the Byzantine agent, the three normal agents update their states mathematically equivalent to the case as if their neighbor graph is the 3-vertex $(1,1)$-reduced graph with the arc set $\{(1,2), (1,3), (2,3)\}$, which is rooted (cf. Lemma \ref{le:convexcombforhighd}). 
In this case, since vertex 1 is the only root and agent 1 does not have any neighbor, it follows the single-agent subgradient algorithm, and thus its state will converge to a minimum point of $f_1(x)$, denoted $x^*$. Since all normal agents will eventually reach a consensus (cf. Lemma~\ref{le:y}), both states of agents 2 and 3 will converge to $x^*$. To guarantee the resilient distributed optimization problem is solvable in this case, there must  hold that $x^*\in\argmin_xf_i(x)$, $i\in\{1,2,3\}$, which implies that the network needs to be 3-redundant. 
It is easy to see that $\kappa_{1,1}(\bbb C)=1$, and thus $n-\kappa_{1,1}(\bbb G)=3$.

Theorem \ref{thm:main} shows that the proposed algorithm achieves full resiliency. 
We next partially characterize the convergence rate of the algorithm.

\begin{theorem}\label{thm:partialrate}
Under Assumptions~\ref{a:interior} and \ref{assum:boundedsubgradient}, if $\bbb G$ is $(\beta,d\beta)$-resilient, the $n$-agent network is $(n-\kappa_{\beta,d\beta}(\bbb G))$-redundant, and $\alpha(t) = 1/\sqrt{T}$ for $T>0$ steps, i.e., $t\in\{0,1,\ldots,T-1\}$,
then there exist a subset of normal agents $\scr S\subset \scr H$ with  $|\scr S| \ge \kappa_{\beta,d\beta}(\bbb G)$, a positive constant $C\ge 1$, and a time subsequence $\scr T\subset \{0,1,\ldots,T-1\}$ with $|\scr T|\ge  T/C$ such that for any $j\in\scr H$ and $x^*\in\scr X^*$, 
\begin{align}
    \sum_{i\in\scr S}f_i\bigg(\frac{\sum_{t\in\scr T}x_j(t)  }{|\scr T|} \bigg) - \sum_{i\in\scr S}f_i(x^*) 
    \le O\Big(\frac{1}{ \sqrt{T}}\Big).\label{eq:jingxuan}
    \end{align}
\end{theorem}

The existing distributed optimization literature (without Byzantine agents) typically characterizes convergence rates by bounding the difference between $\sum_{i\in\scr V}f_i(\frac{1}{T}\sum_{t=0}^{T-1}x_i(t))$ and $\sum_{i\in\scr V}f_i(x^*)$. The above theorem can be viewed as a ``partial'' convergence rate result in that it only reckons a subset $\scr S$ of normal agents and a subsequence $\scr T$ in a finite time horizon. Notwithstanding this, it is worth noting that $\min \sum_{i\in\scr S}f_i(x)$ is equivalent to $\min \sum_{i\in\scr V}f_i(x)$ in the setting here with Byzantine agents (cf. Corollary \ref{coro:optset}) and that $|\scr T|=O(T)$. Therefore, the theorem still to some extent evaluates the convergence rate of the resilient distributed subgradient algorithm under consideration. 
It is well known that the optimal convergence rate of subgradient methods for convex optimization is $O(1/\sqrt{t})$. Whether $f_{\scr H}(\cdot)=\sum_{i\in\scr H} f_i(\cdot)$ converges at this optimal rate or not, has so far eluded us.


Theorem \ref{thm:partialrate} is an immediate consequence of the following proposition. 

\begin{proposition}
\label{thm:jingxuanrate}
Under Assumptions~\ref{a:interior} and \ref{assum:boundedsubgradient}, if $\bbb G$ is $(\beta,d\beta)$-resilient, the $n$-agent network is $(n-\kappa_{\beta,d\beta}(\bbb G))$-redundant, and $\alpha(t) = 1/\sqrt{T}$ for $t\in\{0,1,\ldots,T-1\}$,
then for any integer $b\in[\kappa_{\beta,d\beta}(\bbb G), n-|\scr F|]$, there exist a subset of normal agents $\scr S\subset \scr H$ with  $b\ge |\scr S| \ge \kappa_{\beta,d\beta}(\bbb G)$ and a time subsequence $\scr T\subset \{0,1,\ldots,T-1\}$ with $|\scr T|\ge  T/\sum_{k=\kappa_{\beta,d\beta}(\bbb G)}^{b} \binom{n-|\scr F|}{k}$ such that \eqref{eq:jingxuan} holds for any $j\in\scr H$ and $x^*\in\scr X^*$. 
\end{proposition}






The proposition further quantifies a trade-off between the number of normal agents in $\scr S$ and the length of time subsequence $\scr T$. Roughly speaking, the fewer the normal agents involved in \eqref{eq:jingxuan}, the denser would the time subsequence be. In the special case when $b = \kappa_{\beta,d\beta}(\bbb G)$, the proposition simplifies to the following corollary.

\begin{corollary}\label{thm:rate}
Under Assumptions~\ref{a:interior} and \ref{assum:boundedsubgradient}, if $\bbb G$ is $(\beta,d\beta)$-resilient, the $n$-agent network is $(n-\kappa_{\beta,d\beta}(\bbb G))$-redundant, and $\alpha(t) = 1/\sqrt{T}$ for $t\in\{0,1,\ldots,T-1\}$,
then there exist a subset of normal agents $\scr S\subset \scr H$ with  $|\scr S| = \kappa_{\beta,d\beta}(\bbb G)$ and a time subsequence $\scr T\subset \{0,1,\ldots,T-1\}$ with $|\scr T|\ge  T/\binom{n-|\scr F|}{\kappa_{\beta,d\beta}(\bbb G)}$ such that \eqref{eq:jingxuan} holds for any $j\in\scr H$ and $x^*\in\scr X^*$.
\end{corollary}


\section{Analysis}

This section provides the analysis of the algorithm and proofs of the main results.

\subsection{Algorithm Feasibility}

From Lemma \ref{le:enoughneighbor}, 
$(\beta,d\beta)$-resilient $\bbb G$ guarantees that each agent has at least $(d+1)\beta + 1$ neighbors at each time $t$. Thus, each $\scr A_{ij}$ in the algorithm is always nonempty.

We next show that $y_{ij}(t)$ in \eqref{eq:hull} always exists. To this end, we need the following well-known theorem by Helly.

\begin{lemma}\label{le:helly}
{\rm(Helly's Theorem \cite{Helly1923})} 
Let $\{\scr C_1, \scr C_2\ldots, \scr C_m\}$ be a finite collection of convex sets in $\R^d$ with $m\ge d +1$. If the intersection of every $d+1$ of these sets is nonempty, then the whole collection has a nonempty intersection, i.e.,
$\bigcap_{i=1}^m \scr C_i \neq \emptyset$.
\end{lemma}


\begin{lemma}\label{le:nonempty}
For any $i\in\scr V$ and $j\in\{1,\ldots,a_i\}$, there holds $\bigcap_{\scr S\in\scr B_{ij}} \scr C_{\scr S}(t)\neq\emptyset$.
\end{lemma}

{\bf Proof of Lemma \ref{le:nonempty}:}
From Lemma \ref{le:helly}, it is sufficient to prove that the intersection of every $d+1$ sets in $\scr B_{ij}$ is nonempty. 
Pick any $\scr P \subset \scr B_{ij}$ with $|\scr P|=d+1$.  
For each $\scr S\in\scr P$, 
from the definition of $\scr B_{ij}$, $\scr C_{\scr S}(t)$ is the convex hull of distinct $(d\beta+1)$ points. Since $|\scr P|=d+1$, the intersection $\bigcap_{\scr S\in\scr P}\scr C_{\scr S}(t)$ involves in total $(d+1)(d\beta+1)$ points (with repetition). Recall that all these points are agents' states at time $t$. Thus, each of them can be written as $x_h(t)$ with $h$ being an index in $\scr V$. From the definition of $\scr B_{ij}$, it is easy to see that $h\in\scr A_{ij}$. Note that
\[(d+1)(d\beta+1)-d|\scr A_{ij}| = (d+1)(d\beta+1)-d((d+1)\beta+1)=1.\]
Then, the pigeonhole principle (or Dirichlet's box principle) guarantees that among the total $(d+1)(d\beta+1)$ indices, at least one index in $\scr A_{ij}$, say $k$, repeating at least $(d+1)$ times. Since for each $\scr S\in\scr P$, there is no repetition of indices when computing $\scr C_{\scr S}(t)$, there must exist $(d+1)$ different sets $\scr S_1,\ldots,S_{d+1}\in\scr P$ for which $x_k(t)$ involves the computation of $\scr C_{\scr S_p}(t)$ and thus $x_k(t)\in\scr C_{\scr S_p}(t)$ for all $p\in\{1,\ldots,d+1\}$. 
Since $|\scr P| = d+1$, $\scr P=\{S_1,\ldots,S_{d+1}\}$. It follows that $x_k(t)\in\bigcap_{p=1}^{d+1}\scr C_{\scr S_p}(t)=\bigcap_{\scr S\in\scr P}\scr C_{\scr S}(t)$.
\hfill$\qed$

\subsection{Dynamics of Normal Agents}

To analyze the performance of the algorithm, it is important to understand the dynamics of normal agents and ``decouple'' the influence of Byzantine agents. The following lemma serves this purpose.

\begin{lemma}\label{le:convexcombforhighd}
$v_i(t)$ in \eqref{eq:v_i(t)_ori} can be expressed as a convex combination of $x_i(t)$ and $x_k(t)$, $k\in\scr N_i\cap\scr H$, 
\begin{align*}\label{eq:convexcombhighd}
    v_i(t) = w_{ii}(t)x_i(t) + \sum_{k\in\scr N_i\cap\scr H}w_{ik}(t)x_k(t),
\end{align*}
where $w_{ii}(t)$ and $w_{ik}(t)$ are nonnegative numbers satisfying $w_{ii}(t)+ \sum_{k\in\scr N_i\cap\scr H}w_{ik}(t)=1$, and there exists a positive constant $\eta$ such that for all $i\in\scr H$ and $t$, $w_{ii}(t)\ge \eta$ and among all $w_{ik}(t)$, $k\in\scr N_i\cap\scr H$, at least 
$|\scr N_i\cap\scr H|-d\beta$ of them are bounded below by $\eta$.
\end{lemma}

Since each agent is assumed to have at most $\beta$ Byzantine neighbors, the lemma immediately implies that at least $|\scr N_i|-(d+1)\beta$ among all $w_{ik}(t)$, $k\in\scr N_i\cap\scr H$, are bounded below by $\eta$, which has been reported in \cite[Theorem 1]{resilientconstrained}.
In the special case when $d=1$, the lemma simplifies to Claim~2 in \cite{matrixrepresent}, which directly implies Proposition 5.1 in \cite{sundaram2018distributed}. Thus, the lemma can be regarded as a generalization of \cite[Theorem 1]{resilientconstrained}, \cite[Claim 2]{matrixrepresent}, and \cite[Proposition 5.1]{sundaram2018distributed}.

{\bf Proof of Lemma \ref{le:convexcombforhighd}:}
Recall that there are at most $\beta$ Byzantine neighbors in $\scr A_{ij}$ whose cardinality $|\scr A_{ij}|=(d+1)\beta+1$, and that $\scr B_{ij}$ is the collection of those subsets of $\scr A_{ij}$ whose cardinality is $d\beta+1=|\scr A_{ij}|-\beta$, 
there must exist an index set $\scr P\in \scr B_{ij}$ such that $\scr P\subset\scr N_i\cap\scr H$. For any such index set $\scr P$, since $y_{ij}(t)\in\bigcap_{\scr S\in\scr B_{ij}} \scr C_{\scr S}(t)$, it follows that $y_{ij}(t)\in \scr C_{\scr P}(t)$. 
Let 
$\scr Q_{ij} = \{\scr P:\scr P\in\scr B_{ij}, \scr P\subset\scr N_i\cap\scr H\}$ for each $i\in\scr V$ and $j\in\{1,\ldots,a_i\}$.
From the preceding,  $\scr Q_{ij}$ is always nonempty, and 
for every 
$\scr P\in\scr Q_{ij}$, 
\eq{y_{ij}(t) = \sum_{p\in\scr P}c_p(\scr P)x_p(t),\label{eq:basicconvex}}
where $c_p(\scr P)$ are nonnegative weights satisfying $\sum_{p\in\scr P}c_p(\scr P)=1$. 
It is clear that at least one of $c_p(\scr P)$ is positive and at least $1/|\scr P|=1/(d\beta+1)$.

It is easy to see that $y_{ij}(t)$ can be rewritten as 
\begin{align}\label{eq:expandconvex}
    y_{ij}(t) = \frac{1}{|\scr Q_{ij}|}\sum_{\scr P\in\scr Q_{ij}}\sum_{p\in\scr P}c_p(\scr P)x_p(t).
\end{align}
Our reason for rewriting $y_{ij}(t)$ in this way will be clear shortly.  
Since each $\scr P\subset \scr N_i\cap\scr H$, the above expression is a convex combination of all $x_k(t)$, $k\in\scr N_i\cap\scr H$, allowing some weights being zero.  
Specifically, defining $\scr S_{ijk}=\{\scr P : \scr P\in\scr Q_{ij}, k\in\scr P\}$ for each $\scr Q_{ij}$ and $k\in\scr N_i\cap\scr H$,  
$$y_{ij}(t)=\sum_{k\in\scr N_i\cap\scr H}\Big(\sum_{\scr Q_{ij}}\sum_{\scr P\in\scr S_{ijk}}\frac{c_k(\scr P)}{|\scr Q_{ij}|}\Big)x_k(t).$$
Then, 
\begin{align*}
    v_i(t) &= \frac{1}{1+a_i}\Big(x_i(t)+\sum_{j=1}^{a_i}y_{ij}(t)\Big)\\
    &= \frac{1}{1+a_i}\Big(x_i(t)+\sum_{j=1}^{a_i}\sum_{k\in\scr N_i\cap\scr H}\sum_{\scr Q_{ij}}\sum_{\scr P\in\scr S_{ijk}}\frac{c_k(\scr P)}{|\scr Q_{ij}|}x_k(t)\Big)\\
    & = \frac{x_i(t)}{1+a_i} + \sum_{k\in\scr N_i\cap\scr H}\Big(\sum_{j=1}^{a_i}\sum_{\scr Q_{ij}}\sum_{\scr P\in\scr S_{ijk}}\frac{c_k(\scr P)}{(1+a_i)|\scr Q_{ij}|}\Big)x_k(t)\\
    &=w_{ii}(t)x_i(t) + \sum_{k\in\scr N_i\cap\scr H}w_{ik}(t)x_k(t),
\end{align*}
in which 
\begin{align}
    w_{ii}(t)\dfb\frac{1}{1+a_i},\;\;\;\;\; w_{ik}(t)\dfb \sum_{j=1}^{a_i}\sum_{\scr Q_{ij}}\sum_{\scr P\in\scr S_{ijk}}\frac{c_k(\scr P)}{(1+a_i)|\scr Q_{ij}|},\;\;\;k\in\scr N_i\cap\scr H.\label{eq:convexweights}
\end{align}
It is clear that $w_{ii}(t)>0$ for all $i$ and $t$. 
Since
$$\frac{1}{1+a_i}\sum_{j=1}^{a_i}y_{ij}(t)= \sum_{k\in\scr N_i\cap\scr H}w_{ik}(t)x_k(t),$$
from \eqref{eq:expandconvex}, 
\eq{\sum_{k\in\scr N_i\cap\scr H}w_{ik}(t)x_k(t) = \sum_{j=1}^{a_i}\sum_{\scr P\in\scr Q_{ij}}\sum_{p\in\scr P}\frac{c_p(\scr P)x_p(t)}{(1+a_i)|\scr Q_{ij}|}.\label{eq:finalconvex}}
Note that $\bigcup_{j=1}^{a_i}\bigcup_{\scr P\in\scr B_{ij}}\scr P$ is the collection of all subsets of $\scr N_i$ with cardinality being $d\beta+1$. It follows that $\bigcup_{j=1}^{a_i}\bigcup_{\scr P\in\scr Q_{ij}}\scr P$ is the collection of all subsets of $\scr N_i\cap\scr H$ with cardinality being $d\beta+1$. Thus, \eqref{eq:finalconvex} implies that $\sum_{k\in\scr N_i\cap\scr H}w_{ik}(t)x_k(t)$ is a convex combination of all possible expressions of $y_{ij}(t)$ in terms of \eqref{eq:basicconvex}. Since both $1+a_i$ and $|\scr Q_{ij}|$ are positive, as long as a $c_p(\scr P)$ in \eqref{eq:finalconvex} is positive, $w_{ip}(t)$ is positive with $p\in\scr N_i\cap\scr H$.

We claim that the number of those indices $p\in\scr N_i\cap\scr H$ such that $c_p(\scr P)\ge 1/|\scr P|$ for at least one $\scr P\in \scr Q_{ij}$ is at least $|\scr N_i\cap\scr H|-d\beta$. To prove this, suppose that, to the contrary, the number of such indices is no larger than $|\scr N_i\cap\scr H|-d\beta-1$. 
That is, at least $d\beta+1$ indices in $\scr N_i\cap\scr H$ whose corresponding $c_p(\scr P)<1/|\scr P|$ for all $\scr P\in\scr Q_{ij}$, $j\in\{1,\ldots,a_i\}$. Pick exactly $d\beta+1$ of them and form an index set ${\scr P}_0$. It is clear that $\scr P_0\in\scr Q_{ij}$ for some $j\in\{1,\ldots,a_i\}$. So all $c_p(\scr P_0)$, $p\in\scr P_0$, must be included in the right hand side of \eqref{eq:finalconvex} and strictly less than $1/|\scr P|$. But this is impossible because \eqref{eq:basicconvex} asserts that $\sum_{p\in\scr P_0}c_p(\scr P_0)=1$.

From the preceding, there exist at least $|\scr N_i\cap\scr H|-d\beta$ indices $k\in\scr N_i\cap\scr H$ for which $c_k(\scr P)$ in \eqref{eq:convexweights} is no less than $1/|\scr P|$ for at least one $\scr P\in\scr S_{ijk}$. For each of such $k$, from \eqref{eq:convexweights}, 
\[w_{ik}(t)\ge\frac{1}{|\scr P|(1+a_i)|\scr Q_{ij}|}=\frac{1}{(d\beta+1)(1+a_i)\binom{(d+1)\beta+1}{d\beta+1}}.\]
Set 
\eq{\eta\;\dfb\; \min_{i\in\scr V}\;\frac{1}{(d\beta+1)(1+a_i)\binom{(d+1)\beta+1}{d\beta+1}}.\label{eq:eta}}
Since $a_i\le \binom{n}{(d+1)\beta+1}$ due to \eqref{eq:a_i(t)}, $\eta$ must be positive and independent of $i$ and $t$. 
The statement of the lemma then immediately follows. 
\hfill$\qed$



From \eqref{eq:x_ori} and Lemma \ref{le:convexcombforhighd}, the updates of all normal agents can be written as 
\begin{align}
    v_i(t) &= w_{ii}(t)x_i(t) + \sum_{k\in\scr N_i\cap\scr H}w_{ik}(t)x_k(t),\label{eq:v_i(t)_ana}\\
    x_i(t+1) &= v_i(t)-\alpha(t)g_i(v_i(t)),\label{eq:x_ana}
\end{align}
for all $i\in\scr H$. 

Without loss of generality, we label all normal agents from $1$ to $|\scr H|$ in the sequel. 

Let 
\begin{align*}
    x(t)\dfb \begin{bmatrix}x_1^\top(t)\cr\vdots\cr x_{|\scr H|}^\top(t)\end{bmatrix}, \;\;\;\;\; v(t)\dfb \begin{bmatrix}v_1^\top(t)\cr\vdots\cr v_{|\scr H|}^\top(t)\end{bmatrix}, \;\;\;\;\; g(v(t))\dfb \begin{bmatrix}g_1^\top(v_1(t))\cr\vdots\cr g_{|\scr H|}^\top(v_{|\scr H|}(t))\end{bmatrix}.
\end{align*}
Then, the updates in \eqref{eq:v_i(t)_ana} and \eqref{eq:x_ana} can be written in the form of state equations:
\begin{align}
    v(t)&=W(t)x(t),\label{eq:normalvstate}\\
    x(t+1)&= v(t) - \alpha(t)g(v(t))\label{eq:normalxstate},
\end{align}
where each $W(t)$ is a $|\scr H|\times |\scr H|$ stochastic matrix with positive diagonal entries.\footnote{A square nonnegative matrix is called a stochastic matrix if its row sums all equal one.}

\subsection{Consensus}

We first study the infinite product of stochastic matrices $W(t)$. 

The graph of an $n\times n$ matrix $M$, denoted $\gamma(M)$, is a direct graph with $n$ vertices and a directed edge from vertex $i$ to vertex $j$ whenever the $ji$-th entry of the matrix is nonzero.


\begin{lemma}\label{le:wrooted}
If $\bbb G$ is $(\beta,d\beta)$-resilient, the graph of each $W(t)$ in \eqref{eq:normalvstate} has a rooted spanning subgraph and all the diagonal entries and those off-diagonal entries of $W(t)$ corresponding to the rooted spanning subgraph are uniformly bounded below by a positive number $\eta$ given in \eqref{eq:eta}. 
\end{lemma}

{\bf Proof of Lemma \ref{le:wrooted}:} 
From Lemma \ref{le:convexcombforhighd}, for all $i\in\scr H$, $w_{ii}(t)$ and at least 
$|\scr N_i\cap\scr H|-d\beta$ among all $w_{ik}(t)$, $k\in\scr N_i\cap\scr H$, are positive and uniformly bounded below by $\eta$ given in \eqref{eq:eta}. The expression ``at least 
$|\scr N_i\cap\scr H|-d\beta$ among all $w_{ik}(t)$'' implies that the graph of $W(t)$ can be obtained by removing from each vertex of the subgraph of $\bbb G$ induced by $\scr H$, $\bbb G_{\scr H}$, at most $d\beta$ unspecified incoming edges in $\bbb G_{\scr H}$, which is the $(|\scr F|,d\beta)$-reduced graph of $\bbb G$. Since $|\scr F|\le \beta$, $\bbb G$ must be $(|\scr F|,d\beta)$-resilient. Thus, any $(|\scr F|,d\beta)$-reduced graph of $\bbb G$ is rooted, so is $W(t)$. 
\hfill$\qed$

For any infinite sequence of stochastic matrices with the property in Lemma \ref{le:wrooted}, their product has the following result.

\begin{lemma}\label{le:consensus}
Let $S_1,S_2,\ldots$ be an infinite sequence of $n\times n$ stochastic matrices, each of whose graphs having a rooted spanning subgraph. If all the diagonal entries and those off-diagonal entries of $S_1,S_2,\ldots$ corresponding to the rooted spanning subgraphs are uniformly bounded below by a positive number $p$, then the product $S_k\cdots S_2S_1$ converges to a rank one matrix of the form $\1v^\top$ exponentially fast, where $v$ is a column vector.\footnote{We use $\1$ to denote the vector whose entries all equal to $1$ and the dimension is to be understood from the context.} 
\end{lemma}


To prove the lemma, we  need the concept of the ``composition'' of directed graphs.
The {\em composition} of two directed graphs $\mathbb{G}_p$, $\mathbb{G}_q$
with the same vertex set,   denoted by $\mathbb{G}_q\circ\mathbb{G}_p$,
is the directed graph with the same vertex set and
 arc set defined  so that $(i, j)$ is an arc in the
composition whenever there is a vertex $k$ such that $(i, k)$
 is an arc in $\mathbb{G}_p$ and $(k, j)$
is an arc in $\mathbb{G}_q$.  Since this composition is an associative binary operation, the definition extends unambiguously to any finite sequence of directed graphs with the same vertex set.
Composition and matrix multiplication are closely related.  In particular,  the graph
 of the product of two nonnegative matrices $M_1,M_2\in\R^{n\times n}$
   is equal to the composition of the graphs of the two matrices comprising the product.
    In other words,
   $\gamma(M_2M_1) = \gamma(M_2)\circ\gamma(M_1)$.
If we focus exclusively on graphs with self-arcs at all
vertices, the definition of
composition implies that the arcs of both $\bbb{G}_p$ and $\bbb{G}_q$ are
arcs of $\bbb{G}_q \circ \bbb{G}_p$; the converse is false.

To proceed, for any $n\times n$ nonnegative matrix $S$, let
$$\mu(S)\dfb\max_{i,j}\Big(1-\sum_{k=1}^n\min\{s_{ik},s_{jk}\}\Big).$$
It is easy to see that if $S$ is a substochastic matrix\footnote{A square nonnegative matrix is called substochastic if its row sums are all equal to or less than one.}, $\mu(S)\in[0,1]$. In the case when $S$ is a stochastic matrix, $\mu(S)$ is called the coefficients of ergodicity \cite[page 137]{seneta}. A stochastic matrix $S$ is called a scrambling matrix if and only if $\mu(S)<1$ \cite{hajnal}, whose graph is sometimes called ``neighbor-shared'' \cite{reachingp2}. It is natural to call a vertex $i$ a neighbor of vertex $j$ in a directed graph $\bbb G$ if $(i,j)$ is an arc in $\bbb G$. A directed graph $\bbb G$ is called neighbor-shared if each set of two distinct vertices share a common neighbor. The composition of any set of $n-1$ rooted graphs with $n$ vertices and self-arcs at all vertices, is neighbor shared \cite[Proposition 8]{reachingp1}.
It is easy to check that for any nonnegative square matrix $S$, $\mu(S)<1$ if and only if its graph is neighbor-shared.
Moreover, for any two nonnegative square matrices $S_1$ and $S_2$, if $S_1\ge S_2$,\footnote{For any two real
matrices $A$ and $B$ with the same size, we write $A \ge B$ if $a_{ij}\ge b_{ij}$ for all $i$ and $j$.} then $\mu(S_1)\le \mu(S_2)$.

{\bf Proof of Lemma \ref{le:consensus}:} Since the graph of each $S_t$ is rooted with self-arcs, by Proposition 8 in \cite{reachingp1}, the graph of the product of any finite sequence of $S_t$ matrices of length $n-1$, is neighbor-shared, which implies that the product is a scrambling matrix. 
Thus, letting $V_t=\prod_{\tau=t}^{t+n-2}S_\tau$ for each $t$, each $V_t$ is scrambling and its graph $\gamma(V_t)$ is neighbor-shared. Since the graph of each $S_t$ has a rooted spanning subgraph with self-arcs whose corresponding entries in $S_t$ are bounded below by a positive number $p$. It follows that  the graph of each $V_t$ has a neighbor-shared spanning subgraph $\bbb S$ with self-arcs whose corresponding entries in $V_t$ are bounded below by a positive number $p^{n-1}$. Let $U_t$ be the $n\times n$ matrix whose $ij$th entry is the $ij$th entry of $V_t$ if $(j,i)$ is an arc in $\bbb S$ and zero otherwise. Then, each $U_t$ is a substochastic matrix whose graph is neighbor-shared. Since all positive entries of $U_t$ are bounded below by $p^{n-1}$, 
$\mu(U_t)\le 1-p^{n-1}$.
Since $V_t\ge U_t$,  $\mu(V_t)\le\mu(U_t)\le 1-p^{n-1}$. With this uniform upper bound of all $\mu(V_t)$, the lemma thus is an immediate consequence of Lemma 3 in \cite{hajnal}.
\hfill$\qed$

An infinite sequence of stochastic matrices $\{S(t)\}$ is called ergodic if $\lim_{t\rightarrow\infty}S(t)\cdots S(\tau+1)S(\tau)=\1v^\top(\tau)$ for all $\tau$, where each $v(\tau)$ is a stochastic vector.\footnote{A vector is called a stochastic vector if
its entries are nonnegative and sum to one.}
From Lemmas \ref{le:wrooted} and \ref{le:consensus}, the sequence of stochastic matrices $\{W(t)\}$ is ergodic, and any infinite product of $W(t)$ matrices converges to a rank one matrix 
exponentially fast. Using the same argument as in the proof of Lemma \ref{le:consensus}, 
$$\mu\big(W(t+n-2)\cdots W(t+1)W(t)\big)\le 1-\eta^{n-1}$$
for all $t$, where $\eta$ is given in \eqref{eq:eta}. Following the same argument as pages 610--611 in \cite{reachingp2}, the product $W(t)\cdots W(\tau+1)W(\tau)$ converges to a rank one matrix as $t\rightarrow\infty$ exponentially fast at a rate no slower than
\eq{\lambda\dfb\left(1-\eta^{n-1}\right)^{\frac{1}{n-1}}.\label{eq:lambda}}

Lemmas \ref{le:wrooted} and \ref{le:consensus} also have the following important implication.

\begin{proposition}\label{prop:merelyconsensus}
Under Assumptions \ref{assum:boundedsubgradient} and  \ref{assum:step-size}, if $\bbb G$ is $(\beta,d\beta)$-resilient, all the normal agents will asymptotically reach a consensus.
\end{proposition}

To prove the proposition, we need the following concept.


\begin{definition}\label{def: absolute prob}
Let $\{ S(t) \}$ be a sequence of stochastic matrices. A sequence of stochastic vectors $\{ \pi(t) \}$ is an absolute probability sequence for $\{ S(t) \}$ if
$\pi^\top(t) = \pi^\top(t+1) S(t)$ for all $t\ge0$.
\end{definition}

\vspace{.05in}

This definition was first introduced by Kolmogorov \cite{kolmogorov}. 
It was shown by Blackwell \cite{blackwell} that every sequence of stochastic matrices
has an absolute probability sequence.
In general, a sequence of stochastic matrices may have more than one absolute probability sequence; when the sequence of stochastic matrices is ergodic,
it has a unique absolute probability sequence \cite[Lemma~1]{tacrate}. It is easy to see that when $S(t)$ is a fixed irreducible stochastic matrix $S$, $\pi(t)$ is simply the normalized left eigenvector of $S$ for eigenvalue one, and when $\{S(t)\}$ is an ergodic sequence of doubly stochastic matrices\footnote{A square nonnegative matrix is called a doubly stochastic matrix if its row sums and column sums all equal one.}, $\pi(t)=(1/n)\1$.

From the preceding, the sequence of stochastic matrices $\{W(t)\}$ in \eqref{eq:normalvstate} is ergodic. Thus, $\{W(t)\}$ has a unique absolute probability sequence $\{\pi(t)\}$.  
From Lemma 1 in \cite{tacrate}, \eq{\lim_{t\rightarrow\infty}W(t)\cdots W(\tau+1)W(\tau)=\1\pi^\top(\tau).\label{eq:piexpression}}
Let $\Phi(t,\tau)\dfb W(t)\cdots W(\tau)$ with $t\ge\tau$. Then, there exists a positive constant $c$ such that 
for all $i,j \in \mathcal{H}$ and $t \ge \tau \ge 0$,
\begin{align}\label{eq:shrink}
    \big| \big[\Phi(t,\tau)]_{ij} - \pi_j(\tau) \big|\le c \lambda^{t-\tau},
\end{align}
where $[\cdot]_{ij}$ denotes the $ij$th entry of a matrix and $\lambda$ is given in \eqref{eq:lambda}. Using the same argument as in the proof of Lemma 2 in \cite{nedic}, $c=2$.



To proceed, define
$$y(t)\dfb \pi^\top(t)x(t)=\sum_{i=1}^{|\scr H|} \pi_i(t)x_i(t).$$
From \eqref{eq:x_ana},
\begin{align}
  y(t+1)&= \sum_{i=1}^{|\scr H|} \pi_i(t+1)x_i(t+1)\nonumber\\
  &= \sum_{i=1}^{|\scr H|} \pi_i(t+1)v_i(t)-\alpha(t)\sum_{i=1}^{|\scr H|} \pi_i(t+1)g_i(v_i(t))\nonumber\\
  &=\pi^\top(t+1)W(t)x(t)-\alpha(t)\sum_{i=1}^{|\scr H|} \pi_i(t+1)g_i(v_i(t))\nonumber\\
  &= y(t)-\alpha(t)\sum_{i=1}^{|\scr H|} \pi_i(t+1)g_i(v_i(t)).\label{y}
\end{align}

\begin{lemma}\label{le:y}
Under Assumptions \ref{assum:boundedsubgradient} and  \ref{assum:step-size}, if $\bbb G$ is $(\beta,d\beta)$-resilient,
$\lim_{t\rightarrow\infty} (x_i(t)-y(t)) = 0$ for all $i\in\scr H$.
\end{lemma}

\noindent
{\bf Proof of Lemma \ref{le:y}:}
For all $t> s$,
\begin{eqnarray*}
x_i(t+1) = \sum_{j=1}^{|\scr H|} \left[\Phi(t,s)\right]_{ij}x_j(s)
- \sum_{r=s}^{t-1} \bigg( \sum_{j=1}^{|\scr H|} \left[\Phi(t,r+1)\right]_{ij} \alpha(r)g_j(v_j(r)) \bigg)-\alpha(t)g_i(v_i(t)). \end{eqnarray*}
From \rep{y}, for all $t>s$,
\begin{eqnarray*}
y(t+1) = y(s) - \sum_{r=s}^{t-1}\bigg(\alpha(r)\sum_{j=1}^{|\scr H|} \pi_j(r+1)g_j(v_j(r))\bigg)-\alpha(t)\sum_{i=1}^{|\scr H|} \pi_i(t+1)g_i(v_i(t)).
 \end{eqnarray*}
Set $s=0$.
Then, using \eqref{eq:shrink} and Assumption \ref{assum:boundedsubgradient}, for $t>0$,
\begin{align}
\|x_i(t)-y(t)\| \;\le\; & \sum_{j=1}^{|\scr H|} \big| \left[\Phi(t-1,0)\right]_{ij} - \pi_j(0) \big| \|x_j(0)\| \nonumber\\
& + \sum_{r=0}^{t-2} \sum_{j=1}^{|\scr H|} \big|\left[\Phi(t-1,r+1)\right]_{ij}-\pi_j(r+1)\big| \alpha(r) \|g_j(v_j(r))\| \nonumber\\
& + \alpha(t-1)\|g_i(v_i(t-1))\| + \alpha(t-1)\sum_{i=1}^{|\scr H|} \pi_i(t)\|g_i(v_i(t-1))\| \nonumber\\
\le \;& 2\lambda^t \sum_{j=1}^{|\scr H|} \|x_j(0)\|
+ 2D|\scr H| \sum_{r=0}^{t-2} \lambda^{t-r-2} \alpha(r) + 2D\alpha(t-1) \label{diff}\\
\le \;& 2\lambda^t \sum_{j=1}^{|\scr H|} \|x_j(0)\|
+ 2D|\scr H| \bigg( \sum_{r=0}^{\lfloor \frac{t}{2}\rfloor-1} \lambda^{t-r-2} \alpha(r) + \sum_{r=\lceil \frac{t}{2}\rceil-1}^{t-2} \lambda^{t-r-2} \alpha(r) \bigg)+ 2D\alpha(t-1) \nonumber\\
 \le \;& 2\lambda^t \sum_{j=1}^{|\scr H|} \|x_j(0)\|
+ \frac{2D|\scr H|}{1-\lambda} \Big( \lambda^{\lceil \frac{t}{2}\rceil-1} \alpha(0) + \alpha\big(\lceil \frac{t}{2}\rceil-1\big) \Big)+ 2D\alpha(t-1),\label{eq:stepsize_negative_time}
\end{align}
where $\lfloor\cdot\rfloor$ and $\lceil\cdot\rceil$ denote the floor and ceiling functions, respectively. 
It follows that
$\lim_{t\rightarrow\infty} \|x_i(t)-y(t)\| = 0$.
\hfill$\qed$

The above proof essentially follows the proof of Lemma 8(a) in \cite{nedich2}, generalizing the straight average $y(t)=\frac{1}{|\scr H|}\sum_{i=1}^{|\scr H|} x_i(t)$ to the time-varying weighted average $y(t)=\sum_{i=1}^{|\scr H|} \pi_i(t)x_i(t)$. It can also be proved using the idea of ``input-output consensus stability'' based on a suitably defined semi-norm; see Corollary 1 in \cite{cdc16consensus}.

Proposition \ref{prop:merelyconsensus} is an immediate consequence of Lemma \ref{le:y}. More can be said.

\subsection{Convergence}

\begin{lemma}\label{le:yy}
Under Assumptions \ref{assum:boundedsubgradient} and  \ref{assum:step-size}, if $\bbb G$ is $(\beta,d\beta)$-resilient,
$\sum_{t=0}^{\infty} \alpha(t)\|x_i(t)-y(t)\| < \infty$ for all $i\in\scr H$.
\end{lemma}

{\bf Proof of Lemma \ref{le:yy}:}
From \rep{diff}, for any $t>0$,
\begin{eqnarray*}
\alpha(t)\|x_i(t)-y(t)\| \le 2\alpha(t)\lambda^t \sum_{j=1}^{|\scr H|} \|x_j(0)\|
+ 2D|\scr H| \sum_{r=0}^{t-2} \lambda^{t-r-2} \alpha(t)\alpha(r) + 2D\alpha(t)\alpha(t-1).
\end{eqnarray*}
Since $\alpha(t)\lambda^{t} \le \alpha^2(t) + \lambda^{2t}$ and
$2\alpha(t)\alpha(r) \le \alpha^2(t) + \alpha^2(r)$ for any $t$ and $r$, then for any $t>0$, 
\begin{eqnarray*}
\alpha(t)\|x_i(t)-y(t)\| &\le& 2\alpha^2(t) \sum_{j=1}^{|\scr H|} \|x_j(0)\| + 2 \lambda^{2t} \sum_{j=1}^{|\scr H|} \|x_j(0)\| \cr
&& +\; 2D|\scr H| \alpha^2(t) \sum_{r=0}^{t-2} \lambda^{t-r-2} + 2D|\scr H| \sum_{r=0}^{t-2} \lambda^{t-r-2} \alpha^2(r) + D\big(\alpha^2(t)+\alpha^2(t-1)\big)\cr
&\le & 2 \alpha^2(t) \sum_{j=1}^{|\scr H|} \|x_j(0)\| + 2 \lambda^{2t} \sum_{j=1}^{|\scr H|} \|x_j(0)\| \cr
&& + \;\Big(\frac{2D|\scr H|}{1-\lambda}\Big) \alpha^2(t) + 2D|\scr H| \sum_{r=0}^{t-2} \lambda^{t-r-2} \alpha^2(r)+D\big(\alpha^2(t)+\alpha^2(t-1)\big).
\end{eqnarray*}
Therefore,
\begin{eqnarray*}
\sum_{t=0}^{\infty} \alpha(t)\|x_i(t)-y(t)\| &\le&
2 \bigg(\sum_{t=1}^{\infty}\alpha^2(t)\bigg) \sum_{j=1}^{|\scr H|} \|x_j(0)\| + 2 \bigg(\sum_{t=1}^{\infty}\lambda^{2t}\bigg) \sum_{j=1}^{|\scr H|} \|x_j(0)\| \cr
&& + \Big(\frac{2D|\scr H|}{1-\lambda}\Big) \sum_{t=1}^{\infty}\alpha^2(t) + 2D|\scr H| \sum_{t=1}^{\infty}\sum_{r=0}^{t-2} \lambda^{t-r-2} \alpha^2(r) + 2D\sum_{t=0}^{\infty}\alpha^2(t)\cr
&& + \alpha(0)\|x_i(0)-y(0)\|\cr
&<& \infty,
\end{eqnarray*}
which completes the proof.
\hfill$\qed$

It is known that under some ``standard'' assumptions, if the graphs of an ergodic sequence of stochastic matrices are ``uniformly strongly connected'', all the entries of its unique absolute probability sequence are uniformly bounded below by a positive constant \cite[Theorem 4.8]{touri2012product}. This is not the case here as the sequence of the graphs $W(t)$ matrices may not be ``uniformly strongly connected''.

\begin{proposition}\label{prop:pilowbound}
If $\bbb G$ is $(\beta,d\beta)$-resilient, then for any $t\ge 0$, there exists a subset $\scr S(t)\subset \scr H$ with  $|\scr S(t)|\ge \kappa_{\beta,d\beta}(\bbb G)$ for which all $\pi_i(t)$, $i\in\scr S(t)$, $t\ge 0$ are bounded below by a positive number.
\end{proposition}

To prove the proposition, we need the following concept and results.

We say that a directed graph $\bbb G$ is strongly rooted at vertex $i$ if each other vertex of $\bbb G$ is reachable from vertex $i$ along a directed path of length 1; that is, $\bbb G$ is strongly rooted at $i$ if $i$ is a neighbor of every other vertex in the graph. A directed graph is called strongly rooted if it has at least one vertex at which it is strongly rooted. For a square nonnegative matrix, its graph is strongly rooted at $i$ if and only if its $i$th column is strictly positive. Moreover, for any $n-1$ directed graphs with $n$ vertices and self-arcs which are all rooted at the same vertex $i$, their composition is strongly rooted at $i$ \cite[Proposition 3]{reachingp1}.


{\bf Proof of Proposition \ref{prop:pilowbound}:}
From Lemma \ref{le:wrooted}, the graph of each $W(t)$ is rooted. It is clear that the number of roots in the graph of each $W(t)$ is at least $\kappa_{\beta,d\beta}(\bbb G)$; that is, $\kappa(\gamma(W(t)))\ge \kappa_{\beta,d\beta}(\bbb G)$ for all time $t$. 
Let 
\eq{l \dfb (|\scr H|-\kappa_{\beta,d\beta}(\bbb G)+1)(|\scr H|-2)+1.\label{eq:length}}
We claim that for any finite sequence of $W(t)$ matrices of length $l$, there exists a subset $\scr S\subset \scr H$, depending on the sequence,
with  $|\scr S|\ge \kappa_{\beta,d\beta}(\bbb G)$ such that each $i\in\scr S$ is a root of the graph of some $W(t)$ for at least $|\scr H|-1$ times. To prove the claim, suppose that, to the contrary, such a subset does not exist; that is, at most $\kappa_{\beta,d\beta}(\bbb G)-1$ vertices in $\scr H$ are a root of the graph of some $W(t)$ for at least $|\scr H|-1$ times. For the remaining vertices, they are a root of the graph of some $W(t)$ for at most $|\scr H|-2$ times. Then, the total number of roots of all the graphs of the sequence of $W(t)$ matrices of length $l$ is no large than $(\kappa_{\beta,d\beta}(\bbb G)-1)l+(|\scr H|-\kappa_{\beta,d\beta}(\bbb G)+1)(|\scr H|-2)$. Meanwhile, since $\kappa(\gamma(W(t)))\ge \kappa_{\beta,d\beta}(\bbb G)$ for all time $t$, 
this total number of roots is at least $\kappa_{\beta,d\beta}(\bbb G)l$, which implies that $(\kappa_{\beta,d\beta}(\bbb G)-1)l+(|\scr H|-\kappa_{\beta,d\beta}(\bbb G)+1)(|\scr H|-2)\ge \kappa_{\beta,d\beta}(\bbb G)l$, and thus $l \le (|\scr H|-\kappa_{\beta,d\beta}(\bbb G)+1)(|\scr H|-2)$. But this contradicts \eqref{eq:length}. Therefore, the claim is true. 

The above claim immediately implies that for any time $t\ge 0$ and the corresponding $l$-length stochastic matrix sequence  starting at time $t$, $\{W(t),W(t+1), \ldots,W(t+l-1)\}$, there exists a subset $\scr S(t)\subset \scr H$ with  $|\scr S(t)|\ge \kappa_{\beta,d\beta}(\bbb G)$ such that each $i\in\scr S(t)$ is a root of at least $|\scr H|-1$ graphs among $\gamma(W(t)),\gamma(W(t+1)), \ldots,\gamma(W(t+l-1))$.
From Lemma \ref{le:wrooted}, the graph of each $W(t)$ is rooted with self-arcs. With the two facts of directed graphs with self-arcs that the arcs of each graph in a graph sequence are arcs of their composition, and that any $|\scr H|-1$ graphs with $|\scr H|$ vertices which are all rooted at the same vertex $i$, their composition is strongly rooted at $i$ \cite[Proposition~3]{reachingp1}, 
it follows that for any time $t\ge 0$, the graph of $W(t+l-1)\cdots W(t+1)W(t)$ is strongly rooted at each vertex $i\in\scr S(t)$. Then, each product $W(t+l-1)\cdots W(t+1)W(t)$ has $|\scr S(t)|\ge\kappa_{\beta,d\beta}(\bbb G)$ strictly positive columns whose entries are uniformly bounded below by a positive number $\eta^l$, where $\eta$ is defined in \eqref{eq:eta}. Since each $W(t)$ is a stochastic matrix, it is easy to see that for any $m\ge l$, each product $W(t+m-1)\cdots W(t+1)W(t)$ has $|\scr S(t)|\ge\kappa_{\beta,d\beta}(\bbb G)$ strictly positive columns whose entries are uniformly bounded below $\eta^l$. Then, the statement of the proposition follows from \eqref{eq:piexpression}. 
\hfill$\qed$

We are now in a position to prove Theorem \ref{thm:main} and Proposition \ref{thm:jingxuanrate}.

{\bf Proof of Theorem \ref{thm:main}:}
From \eqref{eq:x_ana} and \eqref{eq:subgradient}, for any $z\in\scr{X}^*$,
\begin{align*}
    \|x_i(t+1)-z\|^2 &= \|v_i(t)-z-\alpha(t)g_i(v_i(t))\|^2 \\
    &\leq \|v_i(t)-z\|^2 + \alpha^2(t)\|g_i(v_i(t))\|^2 - 2\alpha(t)(f_i(v_i(t))-f_i(z)). 
\end{align*}
Since each $W(t)$ is a stochastic matrix and the 2-norm is convex,
$\|v_i(t)-z\|^2 \leq \sum_{j=1}^{|\scr H|} w_{ij}(t)\|x_j(t)-z\|^2$.
Then, with Assumption \ref{assum:boundedsubgradient} and the fact that $\pi_j(t)=\sum_{i=1}^{|\scr H|} \pi_i(t+1)w_{ij}(t)$ for all $j\in\scr H$ and $t$,
\begin{align}
&\;\;\; \sum_{i=1}^{|\scr H|}\pi_i(t+1)\|x_i(t+1)-z\|^2 \nonumber\\
&\leq \sum_{i=1}^{|\scr H|}\pi_i(t+1)\sum_{j=1}^{|\scr H|} w_{ij}(t)\|x_j(t)-z\|^2 +  \alpha^2(t)\sum_{i=1}^{|\scr H|}\pi_i(t+1)\|g_i(v_i(t))\|^2 \nonumber\\
&\;\;\;\;- 2\alpha(t)\sum_{i=1}^{|\scr H|}\pi_i(t+1)(f_i(v_i(t))-f_i(z)) \nonumber \\
&= \sum_{i=1}^{|\scr H|}\pi_i(t)\|x_i(t)-z\|^2 + \alpha^2(t)\sum_{i=1}^{|\scr H|}\pi_i(t+1)\|g_i(v_i(t))\|^2 - 2\alpha(t)\sum_{i=1}^{|\scr H|}\pi_i(t+1)(f_i(v_i(t))-f_i(z))\nonumber \\
&\le \sum_{i=1}^{|\scr H|}\pi_i(t)\|x_i(t)-z\|^2 + \alpha^2(t)D^2 
- 2\alpha(t)\sum_{i=1}^{|\scr H|}\pi_i(t+1)(f_i(v_i(t))-f_i(y(t))) \nonumber\\
&\;\;\;\;- 2\alpha(t)\sum_{i=1}^{|\scr H|}\pi_i(t+1)(f_i(y(t))-f_i(z)). \label{xxx}
\end{align}
Note that
$$|f_i(v_i(t))-f_i(y(t))|\le D\|v_i(t)-y(t)\|
\le D\sum_{j=1}^{|\scr H|} w_{ij}(t)\|x_j(t)-y(t)\|,$$
which implies that
\begin{eqnarray*}
\sum_{i=1}^{|\scr H|}\pi_i(t+1)|f_i(v_i(t))-f_i(y(t))|
&\le& D \sum_{i=1}^{|\scr H|}\pi_i(t+1) \sum_{j=1}^{|\scr H|} w_{ij}(t)\|x_j(t)-y(t)\| \\
&=& D \sum_{i=1}^{|\scr H|}\pi_i(t) \|x_i(t)-y(t)\| \\
&\le& D \sum_{i=1}^{|\scr H|} \|x_i(t)-y(t)\|.
\end{eqnarray*}
From \rep{xxx},
\begin{eqnarray}
\sum_{i=1}^{|\scr H|}\pi_i(t+1)\|x_i(t+1)-z\|^2 
&\leq& \sum_{i=1}^{|\scr H|}\pi_i(t)\|x_i(t)-z\|^2 + \alpha^2(t)D^2 + 2\alpha(t)D \sum_{i=1}^{|\scr H|} \|x_i(t)-y(t)\| \nonumber\\
&&-\; 2\alpha(t)\sum_{i=1}^{|\scr H|}\pi_i(t+1)(f_i(y(t))-f_i(z)).\label{notsure}
\end{eqnarray}
Thus, 
\begin{eqnarray}
2\sum_{t=0}^{\infty}\alpha(t)\sum_{i=1}^{|\scr H|}\pi_i(t+1)(f_i(y(t))-f_i(z)) &\leq& \sum_{i=1}^{|\scr H|}\pi_i(0)\|x_i(0)-z\|^2 + D^2\sum_{t=0}^{\infty}\alpha^2(t) \nonumber\\
&&+\; 2D\sum_{t=0}^{\infty}\alpha(t) \sum_{i=1}^{|\scr H|} \|x_i(t)-y(t)\| \nonumber\\
&<& \infty \nonumber 
\end{eqnarray}
because of Assumption \ref{assum:step-size} and Lemma \ref{le:yy}.

From Corollary \ref{le:localoptimal}, 
$z\in \argmin_{x}f_i(x)$ for all $i\in\scr H$. Then, $f_i(y(t))-f_i(z)\ge 0$ for all $i\in\scr H$ and $t\ge 0$. 
From Proposition \ref{prop:pilowbound} and its proof, for any $t\ge 0$, there exists a subset $\scr S(t)\subset \scr H$ with  $|\scr S(t)|\ge \kappa_{\beta,d\beta}(\bbb G)$ for which all $\pi_i(t)$, $i\in\scr S(t)$, $t\ge 0$ are uniformly bounded below by $\eta^l>0$.  
It follows that
\begin{equation}\label{eq:eta^landpi}
\eta^l\sum_{t=0}^{\infty}\alpha(t)\sum_{i\in\scr S(t)}(f_i(y(t))-f_i(z))\le 
\sum_{t=0}^{\infty}\alpha(t)\sum_{i=1}^{|\scr H|}\pi_i(t+1)(f_i(y(t))-f_i(z)) < \infty.
\end{equation}
Since $\sum_{t=0}^\infty \alpha(t)=\infty$ and $f_i(y(t))-f_i(z)\ge 0$ for all $i\in\scr H$ and $t\ge 0$,
\[\liminf_{t\rightarrow\infty}\sum_{i\in\scr S(t)}\big(f_i(y(t))-f_i(z)\big)=0.\]
which implies that\footnote{We use ${\rm dist}(x,\scr S)\dfb\inf \{\|x-y\|:y\in\scr S\}$ to denote the Euclidean distance between a point $x$ and a set $\scr S$ in $\R^d$.} \[\liminf_{t\rightarrow\infty}\; {\rm dist} \Big(y(t), \; \argmin_{x} \sum_{i\in\scr S(t)} f_i(x)\Big)=0.\]
From Corollary~\ref{coro:optset}, $\argmin_{x} \sum_{i\in\scr S(t)} f_i(x) = \scr X^*$ for all time $t\ge 0$. Then, 
\begin{align*}
    \argmin_x \Big(\liminf_{t\rightarrow\infty} \sum_{i\in\scr S(t)} f_i(x) \Big)
    = \argmin_x f(x),
\end{align*}
which implies that $\liminf_{t\rightarrow\infty}{\rm dist}(y(t),\scr X^*)=0$, and thus
\begin{align}\label{eq:f_no_liminf}
    \liminf_{t\rightarrow\infty} \big(f(y(t))-f(z)\big)=0.
\end{align}

We next show that all the sequences $\{x_i(t)\}$, $i\in\scr H$, converge to the
same optimal point.

From \eqref{notsure}, rearranging the terms and fixing an arbitrary period from time $t_1$ to $t_2$ with $t_1<t_2$,
\begin{eqnarray}
\sum_{i=1}^{|\scr H|}\pi_i(t_2+1)\|x_i(t_2+1)-z\|^2 
&\leq& \sum_{i=1}^{|\scr H|}\pi_i(t_1)\|x_i(t_1)-z\|^2 \nonumber\\
&&+\; D^2\sum_{t=t_1}^{t_2} \alpha^2(t) + 2D\sum_{t=t_1}^{t_2} \alpha(t) \sum_{i=1}^{|\scr H|} \|x_i(t)-y(t)\|. \nonumber
\end{eqnarray}
From Assumption \ref{assum:boundedsubgradient} and Lemma \ref{le:yy}, 
$$\limsup_{\tau_2\rightarrow\infty}\; \sum_{i=1}^{|\scr H|}\pi_i(\tau_2+1)\|x_i(\tau_2+1)-z\|^2 \le \liminf_{\tau_1\rightarrow\infty}\; \sum_{i=1}^{|\scr H|}\pi_i(\tau_1)\|x_i(\tau_1)-z\|^2.$$
Thus, the sequence $\{\sum_{i=1}^{|\scr H|}\pi_i(t)\|x_i(t)-z\|^2\}$ is convergent
for each $z\in\scr X^*$. 
From Proposition~\ref{prop:pilowbound} and its proof, 
$$\eta^l\sum_{i\in\scr S(t)}\|x_i(t)-z\|^2\le \sum_{i\in\scr S(t)}\pi_i(t)\|x_i(t)-z\|^2\le \sum_{i=1}^{|\scr H|}\pi_i(t)\|x_i(t)-z\|^2,$$ 
which implies that the sequence $\{\sum_{i\in\scr S}\|x_i(t)-z\|^2\}$ is bounded, so is  each sequence $\{x_i(t)\}$, $i\in\scr S(t)$. From Proposition~\ref{prop:merelyconsensus}, 
all the sequences $\{x_i(t)\}$, $i\in\scr H$, are bounded. 
From Lemma \ref{le:y}, the sequence $\{y(t)\}$ is bounded, and with each $\pi(t)$ being a stochastic vector, the sequence $\{\|y(t)-z\|^2\}$ is convergent for each $z\in\scr X^*$. 
Since $y(t)$ is bounded, from \eqref{eq:f_no_liminf}, there exists a subsequence of $\{y(t)\}$ converging to a point $x^*\in\scr X^*$. Since $\{\|y(t)-x^*\|^2\}$ is convergent, the sequence $\{y(t)\}$ converges to $x^*$. From Lemma \ref{le:y}, all the sequences $\{x_i(t)\}$, $i\in\scr H$ converge to the same optimal point $x^*$.
\hfill$\qed$



{\bf Proof of Proposition~\ref{thm:jingxuanrate}:}
From Proposition \ref{prop:pilowbound}, for any $t\ge 0$, there exists a subset $\scr S(t)\subset \scr H$ with  $|\scr S(t)|\ge \kappa_{\beta,d\beta}(\bbb G)$ for which all $\pi_i(t)$, $i\in\scr S(t)$, $t\ge 0$ are bounded below by a positive number.
It is clear that $|\scr S(t)|\le |\scr H|=n-|\scr F|$ for all $t\ge 0$. 
Fix any integer $b\in[\kappa_{\beta,d\beta}(\bbb G), n-|\scr F|]$. 
Define 
\begin{align*}
\scr S = \argmax_{\scr K\subset \scr H,\; b\ge |\scr K|\ge \kappa_{\beta,d\beta}(\bbb G)} \big|\big\{t: t\in\{0,1,\ldots,T-1\}, \;\scr K\subset\scr S(t)\big\}\big|,
\end{align*}
which represents a subset of normal agents, with an appropriate cardinality, appearing the most times within $T$ time steps. It is easy to see that $\scr S$ is well-defined and nonempty. 

We next count the number of times such an $\scr S$ set appears within $T$ time steps. To this end, define \begin{align*}
\scr T = \big\{t:t\in\{0,1,\ldots,T-1\},\; \scr S\subset \scr S(t)\big\},
\end{align*}
which represents this number. Consider the set $\scr M=\{\scr K\subset \scr H: b\ge |\scr K|\ge \kappa_{\beta,d\beta}(\bbb G)\}$ whose cardinality is 
$$C=\sum_{k=\kappa_{\beta,d\beta}(\bbb G)}^{b} \binom{n-|\scr F|}{k}.$$
Note that $\scr S\in\scr M$ and for each $t\ge 0$, there exists at least one set in $\scr M$ being a subset of $\scr S(t)$. 
From the pigeonhole principle, at least one set in $\scr M$ appears $\lceil T/C \rceil$ times within $T$ time steps.
Since $\scr S$ is a set in $\scr M$ that appears most often during the total $T$ time steps, $|\scr T|\ge T/C$.


To proceed, similar to \eqref{eq:eta^landpi}, for any $z\in\scr{X}^*$,
\begin{align}\label{eq:tildeS}
\eta^l\sum_{t=0}^{T-1}\alpha(t)\sum_{i\in{\scr S}(t)}(f_i(y(t))-f_i(z))\le 
\sum_{t=0}^{T-1}\alpha(t)\sum_{i\in\scr H}\pi_i(t+1)(f_i(y(t))-f_i(z)) < \infty.
\end{align}

Substituting $\alpha(t)=1/\sqrt{T}$ to \eqref{notsure} and \eqref{eq:tildeS}, it follows that for any $z\in\scr X^*$, 
\begin{align}\label{eq:fiy-fiz}
     2\eta^l\sum_{t=0}^{T-1}\sum_{i\in{\scr S}(t)}\frac{f_i(y(t))-f_i(z)}{\sqrt{T}}&\le\sum_{i\in\scr H}\pi_i(0)\|x_i(0)-z\|^2 + D^2+ 2D\sum_{t=0}^{T-1} \sum_{i\in\scr H}\frac{\|x_i(t)-y(t)\|}{\sqrt{T}}.
\end{align}
Since all $f_i(x)$, $i\in\scr H$ are convex functions, 
 \begin{align}
     \sum_{i\in\scr S}f_i\bigg(\frac{\sum_{t\in\scr T}y(t)}{|\scr T|}\bigg)-\sum_{i\in\scr S}f_i(z)
     \le \sum_{i\in\scr S}\sum_{t\in\scr T}\frac{f_i(y(t))-f_i(z)}{|\scr T|}
     =\sum_{t\in\scr T}\sum_{i\in\scr S}\frac{f_i(y(t))-f_i(z)}{|\scr T|}.\label{eq:middleconvex}
\end{align}
From Corollary \ref{le:localoptimal}, 
$z\in \argmin_{x}f_i(x)$ for all $i\in\scr H$. Then, $f_i(y(t))-f_i(z)\ge 0$ for all $i\in\scr H$ and $t\ge 0$. Since $\scr S\subset\scr S(t)\subset \scr H$ for all $t\ge 0$, it follows from \eqref{eq:middleconvex}, \eqref{eq:fiy-fiz}, and $|\scr T|\ge T/C$ that 
\begin{align}
     &\sum_{i\in\scr S}f_i\bigg(\frac{\sum_{t\in\scr T}y(t)}{|\scr T|}\bigg)-\sum_{i\in\scr S}f_i(z)
     \le \sum_{t=0}^{T-1}\sum_{i\in\scr S(t)}\frac{f_i(y(t))-f_i(z)}{|\scr T|} \nonumber\\
     \le\;& \frac{C}{2\eta^{l}}\bigg(\frac{\sum_{i\in\scr H}\pi_i(0)\|x_i(0)-z\|^2 + D^2}{\sqrt{T}}+\frac{2D\sum_{t=0}^{T-1} \sum_{i\in\scr H} \|x_i(t)-y(t)\|}{T}\bigg).\label{eq:foryixuanapproach}
 \end{align}
Note that for all $i\in\scr H$, there holds 
$$\|x_i(0)-y(0)\|=\|x_i(0)-\sum_{j\in\scr H}\pi_i^\top(0)x_i(0)\|\le \sum_{j\in\scr H}\|x_j(0)\|.$$
With this simple fact and $\alpha(t)=1/\sqrt{T}$, it follows from \eqref{eq:stepsize_negative_time} that for all $i\in\scr H$,
\begin{align}
\sum_{t=0}^{T-1} \|x_i(t)-y(t)\|\le\;&
\sum_{t=1}^{T-1} 2\lambda^t \sum_{j\in\scr H} \|x_j(0)\|
+ \frac{2D|\scr H|}{1-\lambda} \sum_{t=1}^{T-1} \Big( \lambda^{\lceil \frac{t}{2}\rceil-1} \alpha(0) + \alpha\big(\lceil \frac{t}{2}\rceil-1\big) \Big)+ 2D\sum_{t=1}^{T-1} \alpha(t-1) \nonumber\\
&+\|x_i(0)-y_i(0)\|\nonumber\\
\le\;& \sum_{t=0}^{T-1} 2\lambda^t \sum_{j\in\scr H} \|x_j(0)\|
+ \frac{2D|\scr H|}{1-\lambda} \sum_{t=1}^{T-1} \Big( \lambda^{\lceil \frac{t}{2}\rceil-1} \alpha(0) + \alpha\big(\lceil \frac{t}{2}\rceil-1\big) \Big)+ 2D\sum_{t=1}^{T-1} \alpha(t-1) \nonumber\\
\le\;& \frac{2}{1-\lambda} \sum_{j\in\scr H}\|x_j(0)\|
+ \frac{2D|\scr H|}{1-\lambda} \Big( \frac{2}{\lambda(1-\lambda)\sqrt{T}} + \sqrt{T} \Big)+ 2D\sqrt{T}, \label{eq:lemma9exten}
\end{align}
which implies that
 \begin{align}
  \sum_{i\in\scr S}f_i\bigg(\frac{\sum_{t\in\scr T}y(t)}{|\scr T|}\bigg)-\sum_{i\in\scr S}f_i(z)
  \le\;& \frac{C}{2\eta^{l}}\bigg(\frac{\sum_{i\in\scr H}\pi_i(0)\|x_i(0)-z\|^2 + D^2}{\sqrt{T}}+
  \frac{4D |\scr H| }{(1-\lambda)T} \sum_{j=1}^{|\scr H|} \|x_j(0)\| \nonumber\\
  &\;\;\;\;\;\;\;\; + \frac{8|\scr H|^2D^2}{\lambda(1-\lambda)^2 T\sqrt{T}} + \frac{4|\scr H|^2D^2}{(1-\lambda)\sqrt{T}} + \frac{4D^2 |\scr H|  }{\sqrt{T}}
  \bigg). \label{qqqq}
 \end{align}
From \eqref{eq:subgradient} and Assumption~\ref{assum:boundedsubgradient},
for all $j\in\scr H$, it holds that
\begin{align*}
    &\sum_{i\in\scr S}f_i\bigg(\frac{\sum_{t\in\scr T}x_j(t)}{|\scr T|}\bigg)-\sum_{i\in\scr S}f_i\bigg(\frac{\sum_{t\in\scr T}y(t)}{|\scr T|}\bigg)\\
    \le\;& D\sum_{i\in\scr S} \frac{\sum_{t\in\scr T} \|x_j(t)-y(t)\|}{|\scr T|}
    \le D \sum_{i\in\scr S} \frac{\sum_{t=0}^{T-1} \|x_j(t)-y(t)\|}{|\scr T|}\\
    \le\;& DC|\scr S|\bigg(\frac{2}{(1-\lambda)T} \sum_{j\in\scr H} \|x_j(0)\|
    + \frac{4D|\scr H|}{\lambda(1-\lambda)^2T\sqrt{T}}  
    + \frac{2D|\scr H|}{(1-\lambda)\sqrt{T}}  
    + \frac{2D}{\sqrt{T}}\bigg),
\end{align*}
which, with \eqref{qqqq}, implies that
\begin{align}
  & \sum_{i\in\scr S}f_i\bigg(\frac{\sum_{t\in\scr T}x_j(t)}{|\scr T|}\bigg)-\sum_{i\in\scr S}f_i(z)  
  \nonumber\\
  \le\;& \frac{C}{2\eta^{l}}\bigg(\frac{\sum_{i\in\scr H}\pi_i(0)\|x_i(0)-z\|^2 + D^2}{\sqrt{T}}+
  \frac{4D |\scr H| }{(1-\lambda)T} \sum_{j\in\scr H} \|x_j(0)\| + \frac{8|\scr H|^2D^2}{\lambda(1-\lambda)^2 T\sqrt{T}} + \frac{4|\scr H|^2D^2}{(1-\lambda)\sqrt{T}}  \nonumber\\
  & + \frac{4D^2 |\scr H|  }{\sqrt{T}}
  \bigg) + DC|\scr S| \bigg(\frac{2}{(1-\lambda)T} \sum_{j\in\scr H} \|x_j(0)\|
    + \frac{4D|\scr H|}{\lambda(1-\lambda)^2T\sqrt{T}}  
    + \frac{2D|\scr H|}{(1-\lambda)\sqrt{T}}  
    + \frac{2D}{\sqrt{T}}\bigg)\nonumber\\
  =\;& \frac{C}{2\eta^{l}}\bigg(\frac{\sum_{i\in\scr H}\pi_i(0)\|x_i(0)-z\|^2 + D^2}{\sqrt{T}}+
  \frac{4D \big(|\scr H|+ |\scr S|\eta^l\big) }{(1-\lambda)T} \sum_{j\in\scr H} \|x_j(0)\| 
  + \frac{8|\scr H|D^2\big(|\scr H|+ |\scr S|\eta^l\big)}{\lambda(1-\lambda)^2 T\sqrt{T}} 
   \nonumber\\
  & + \frac{4|\scr H|D^2\big(|\scr H|+ |\scr S|\eta^l\big)}{(1-\lambda)\sqrt{T}}  + \frac{4D^2 \big(|\scr H|+ |\scr S|\eta^l\big)  }{\sqrt{T}} 
    \bigg) \label{aaaa}\\
  =\;& O\Big(\frac{1}{ \sqrt{T}}\Big).\nonumber
 \end{align}

There is an alternative way to prove this, as follows.   

From \rep{notsure}, for any $j\in\scr H$, 
\begin{align*}
    &\;\;\;\; \sum_{i=1}^{|\scr H|}\pi_i(t+1)\|x_i(t+1)-z\|^2 \\
    &\leq \sum_{i=1}^{|\scr H|}\pi_i(t)\|x_i(t)-z\|^2 + \alpha^2(t)D^2 + 2\alpha(t)D \sum_{i=1}^{|\scr H|} \|x_i(t)-y(t)\| \\
    &\;\;\;\; - 2\alpha(t)\sum_{i=1}^{|\scr H|}\pi_i(t+1)(f_i(x_j(t))-f_i(z)) - 2\alpha(t)\sum_{i=1}^{|\scr H|}\pi_i(t+1)(f_i(y(t))-f_i(x_j(t)))\\
    &\leq \sum_{i=1}^{|\scr H|}\pi_i(t)\|x_i(t)-z\|^2 + \alpha^2(t)D^2 + 2\alpha(t)D \sum_{i=1}^{|\scr H|} \|x_i(t)-y(t)\| \\
    &\;\;\;\; - 2\alpha(t)\sum_{i=1}^{|\scr H|}\pi_i(t+1)(f_i(x_j(t))-f_i(z)) 
     + 2D\alpha(t)\sum_{i=1}^{|\scr H|}\pi_i(t+1)\| y(t)-x_j(t) \| \\
    & = \sum_{i=1}^{|\scr H|}\pi_i(t)\|x_i(t)-z\|^2 + \alpha^2(t)D^2 + 2\alpha(t)D \sum_{i=1}^{|\scr H|} \|x_i(t)-y(t)\| \\
    &\;\;\;\; + 2D\alpha(t)\| y(t)-x_j(t) \| - 2\alpha(t)\sum_{i=1}^{|\scr H|}\pi_i(t+1)(f_i(x_j(t))-f_i(z)),
\end{align*}
which implies that
\begin{align*}
    2\alpha(t)\sum_{i=1}^{|\scr H|}\pi_i(t+1)(f_i(x_j(t))-f_i(z))
    & \le - \sum_{i=1}^{|\scr H|}\pi_i(t+1)\|x_i(t+1)-z\|^2 + \sum_{i=1}^{|\scr H|}\pi_i(t)\|x_i(t)-z\|^2  \\
    &\;\;\;\; + \alpha^2(t)D^2 + 2\alpha(t)D \sum_{i=1}^{|\scr H|} \|x_i(t)-y(t)\| \\
    &\;\;\;\; + 2D\alpha(t)\| y(t)-x_j(t) \| .
\end{align*}
Thus,
\begin{align}
    &\;\;\;\;\; 2\sum_{t=0}^{T-1}\alpha(t)\sum_{i=1}^{|\scr H|}\pi_i(t+1)(f_i(x_j(t))-f_i(z)) \nonumber\\
    &\leq \sum_{i=1}^{|\scr H|}\pi_i(0)\|x_i(0)-z\|^2 + D^2\sum_{t=0}^{T-1}\alpha^2(t) + 2D\sum_{t=0}^{T-1}\alpha(t) \sum_{i=1}^{|\scr H|} \|x_i(t)-y(t)\|  + 2D\sum_{t=0}^{T-1}\alpha(t) \|x_j(t)-y(t)\|. \label{eq:fxj}
\end{align}
Similar to \eqref{eq:eta^landpi}, for any $z\in\scr X^*$,
\begin{align}
    \eta^l\sum_{t=0}^{T-1} \alpha(t)\sum_{i\in\scr S(t)}(f_i(x_j(t))-f_i(z))
    \le \sum_{t=0}^{T-1} \alpha(t)\sum_{i=1}^{|\scr H|}\pi_i(t+1)(f_i(x_j(t))-f_i(z)). \label{eq:eta_fxj}
\end{align}
Substituting $\alpha(t)=1/\sqrt{T}$ to \eqref{eq:fxj} and \eqref{eq:eta_fxj}, it holds that for any $z\in\scr X^*$,  
\begin{align*}
    & \sum_{t=0}^{T-1} \sum_{i\in\scr S(t)}(f_i(x_j(t))-f_i(z))
    \le \eta^{-l} \sum_{t=0}^{T-1} \sum_{i=1}^{|\scr H|}\pi_i(t+1)(f_i(x_j(t))-f_i(z))\\
    \le\;& \frac{\sqrt{T}}{2\eta^{l}} \bigg( \sum_{i=1}^{|\scr H|}\pi_i(0)\|x_i(0)-z\|^2 + D^2 + \frac{2D}{\sqrt{T} }\sum_{t=0}^{T-1} \Big( \sum_{i=1}^{|\scr H|} \|x_i(t)-y(t)\|  + \|x_j(t)-y(t)\| \Big) \bigg).
\end{align*}
Using the same argument as in deriving \eqref{eq:foryixuanapproach}, 
\begin{align*}
     & \sum_{i\in\scr S}f_i\bigg(\frac{\sum_{t\in\scr T}x_j(t)}{|\scr T|}\bigg)-\sum_{i\in\scr S}f_i(z)\\
     \le\;& \frac{C}{2\eta^{l}}\bigg(\frac{\sum_{i=1}^{|\scr H|}\pi_i(0)\|x_i(0)-z\|^2 + D^2}{\sqrt{T}}+\frac{2D\sum_{t=0}^{T-1} \big( \sum_{i=1}^{|\scr H|} \|x_i(t)-y(t)\|  + \|x_j(t)-y(t)\| \big)}{T}\bigg), 
 \end{align*}
which, with \eqref{eq:lemma9exten}, yields 
\begin{align}
     & \sum_{i\in\scr S}f_i\bigg(\frac{\sum_{t\in \scr T} x_j(t)}{T}\bigg)-\sum_{i\in\scr S}f_i(z)\nonumber\\
     \le\;& \frac{C}{2\eta^{l}}\bigg(\frac{\sum_{i=1}^{|\scr H|}\pi_i(0)\|x_i(0)-z\|^2 + D^2}{\sqrt{T}}
     +\frac{4D(|\scr H|+1)}{(1-\lambda)T} \sum_{j=1}^{|\scr H|} \|x_j(0)\|
     + \frac{8D^2|\scr H| (|\scr H|+1)}{\lambda(1-\lambda)^2T\sqrt{T}}  \nonumber\\
     &\;\;\;\;\;\;\;\; + \frac{4D^2|\scr H| (|\scr H|+1)}{(1-\lambda)\sqrt{T}} 
     + \frac{4D^2(|\scr H|+1)}{\sqrt{T}} 
     \bigg)\label{zzzz}\\
     =\;& O\Big(\frac{1}{ \sqrt{T}}\Big).\nonumber
 \end{align}
This completes the proof.
\hfill$\qed$


The above proof provides two bounds with the same order. It is not hard to check that the first bound \eqref{aaaa} is better than the second one \eqref{zzzz} if and only if $\eta^l|\scr S|<1$, and they are equal if and only if $\eta^l|\scr S|=1$.

\section{Concluding Remarks}

This paper has proposed a distributed subgradient algorithm which achieves full resilience in the presence of Byzantine agents, with appropriate redundancy in both graph connectivity and objective functions. The algorithm and convergence results can be easily extended to time-varying neighbor graphs, provided that the neighbor graph is $(\beta,d\beta)$-resilient all the time. One immediate next step is to relax Assumption \ref{a:interior}, possibly appealing to gradient descent for differentiable convex functions. The concepts and tools developed in the paper are expected to be applicable to other consensus-based distributed optimization and computation problems. 

Although the algorithm theoretically works for multi-dimensional convex optimization, it has the following limitations which preclude its applicability to high-dimensional optimization. First, from Lemma \ref{le:enoughneighbor}, the algorithm implicitly requires that each agent have at least $(d+1)\beta+1$ neighbors, which is impossible for high dimensions. Second, picking a point in the intersection of multiple convex hulls (cf. step \eqref{eq:hull} in the algorithm) can be computationally expensive in high dimensions, although the issue has been attenuated in \cite[Algorithm 2]{resilientconstrained} and \cite[Section 5.1]{yan2022resilient}. Last, building $(\beta,d\beta)$-resilient graphs is not an easy job, especially when $d$ or $\beta$ is large. Another practical issue of the algorithm, independent of dimensions, is how to measure and establish objective function redundancy. Studies of $(r,s)$-resilient graphs and $k$-redundant multi-agent networks are of independent interest.

Considering that nowadays distributed optimization algorithms in machine learning are frequently high-dimensional,  there is ample motivation to design fully resilient high-dimensional distributed optimization algorithms. A future direction of this paper aims to tackle this challenging problem by combining the proposed algorithm with communication-efficient schemes in which each agent can transmit only low-dimensional signals. Possible approaches include entry-wise or block-wise updating \cite{cdc20,block}, limited information fusion \cite{cdc22}, and dimension-independent filtering \cite{gupta2021byzantine,acc22}.



\bibliographystyle{unsrt}
\bibliography{push,resilience}
\end{document}